\documentclass[12pt, a4paper]{article}
\usepackage[english]{babel}
\usepackage[utf8]{inputenc}

\usepackage{amsmath}
\usepackage{amssymb}
\usepackage{amsfonts}
\usepackage{amsbsy}

\usepackage{bm}

\usepackage{graphicx}

\usepackage{newlfont}
\usepackage{float}
\usepackage[paperwidth=192mm,
			paperheight=262mm,
			vmargin={19mm,19mm},
			hmargin={13.7mm,13.7mm},
			headsep=12pt,
			footskip=12pt]{geometry}

\usepackage{hyperref}

\usepackage{siunitx}

\usepackage{subcaption}

\usepackage{caption}
\usepackage{booktabs}
\usepackage{tabularx}

\DeclareMathOperator{\grad}{\nabla}
\DeclareMathOperator{\dive}{\nabla\cdot}

\begin{document}

\title{An implicit DG solver for incompressible two-phase flows with an artificial compressibility formulation}

\author{Giuseppe Orlando$^{(1),}$\thanks{\textbf{Present address}: CMAP, CNRS, \'{E}cole polytechnique, Institut Polytechnique de Paris, Route de Saclay, 91120 Palaiseau, France, {\tt giuseppe.orlando@polytechnique.edu}}}

\date{}

\maketitle

\begin{center}
{
\small
$^{(1)}$  
MOX - Dipartimento di Matematica, Politecnico di Milano \\
Piazza Leonardo da Vinci 32, 20133 Milano, Italy \\
{\tt giuseppe.orlando@polimi.it}
}
\end{center}


\vspace*{0.5cm}

\pagebreak

\abstract{We propose an implicit Discontinuous Galerkin (DG) discretization for incompressible two-phase flows using an artificial compressibility formulation. The conservative level set (CLS) method is employed in combination with a reinitialization procedure to capture the moving interface. A projection method based on the L-stable TR-BDF2 method is adopted for the time discretization of the Navier-Stokes equations and of the level set method. Adaptive Mesh Refinement (AMR) is employed to enhance the resolution in correspondence of the interface between the two fluids. The effectiveness of the proposed approach is shown in a number of classical benchmarks. A specific analysis on the influence of different choices of the mixture viscosity is also carried out.}

\pagebreak

\section{Introduction}
\label{sec:intro} \indent

Two-phase flows are common in many engineering and industrial applications. An evolving interface delimits the bulk regions of the single phases. Many techniques
have been developed over the years to capture the motion of the interface. Two classes of methods are commonly used to locate the interface: interface-tracking and interface-capturing. Interface-tracking schemes employ either Arbitrary Lagrangian–Eulerian (ALE) methods on a mesh that deforms with the interface \cite{dettmer:2003, hughes:1981} or marker and cell methods \cite{rudman:1997}. Interface-capturing techniques are instead based on fixed spatial grids with an interface function which captures the interface. A full survey on interface-capturing methods goes beyond the scope of this work and we refer e.g. to \cite{mirjalili:2017} for a review of these techniques. Interface capturing methods include the level set (LS) method \cite{osher:2005, osher:1988}, which represents the interface as an iso-surface of the so-called level set function. Classically, the level set function is defined as the signed distance function. However, this choice leads to non conservative methods. A number of approaches have been developed to overcome this issue; in this work, we employ the conservative level set (CLS) method, originally proposed in \cite{olsson:2005, olsson:2007}, and briefly summarized in Section \ref{ssec:clsm}. CLS includes a reinitialization equation to maintain the shape of the level set, which will be also discussed in Section \ref{ssec:clsm}.

Changing fluid properties, such as density and viscosity, and surface tension at the interface lead to discontinuities that make the discretization of the Navier-Stokes equations particularly challenging. The Discontinuous Galerkin (DG) method has been widely employed in the field of Computational Fluid Dynamics (CFD), see e.g. \cite{bassi:1997, fehn:2017, karniadakis:2005}, and is a natural candidate for the discretization of the governing equations of two-phase flows. Several approaches have been proposed in the literature combining the DG method and the level set method, see among many others \cite{gao:2018, grooss:2006, owkes:2013, pochet:2013}. In this paper, we propose an extension of the solver for single-phase incompressible Navier-Stokes equations with an artificial compressibility formulation presented in \cite{orlando:2023b, orlando:2022}, so as to overcome well know issues of projection methods. The time discretization is therefore based on the TR-BDF2 scheme \cite{bank:1985, bonaventura:2017, hosea:1996, orlando:2022}, which is a second order two-stage method. A brief review of the TR-BDF2 method will be given in Section \ref{sec:time_disc}, whereas we refer to \cite{bonaventura:2017, hosea:1996} for a detailed analysis of the scheme. The solver is implemented in the framework of the open source numerical library \texttt{deal.II} \cite{arndt:2022}, which supports native non-conforming \(h-\)adaptation. We will exploit these capabilities to enhance the resolution in the regions close to the interface between the two fluids.

The paper is structured as follows: the model equations and their non-dimensional formulation are reviewed in Section \ref{sec:modeleq}. The time discretization approach is outlined and discussed in Section \ref{sec:time_disc}. The spatial discretization is presented in Section \ref{sec:space_disc}. The application of the proposed method to a number of significant benchmarks is reported in Section \ref{sec:tests}. Here, we also analyze the impact of different possible choices for the mixture viscosity when the interface undergoes large deformations. Finally, some conclusions and perspectives for future work are presented in Section \ref{sec:conclu}.
 
%%%%%%%%%%%%%%%%%%%%%%%%%%%%%%%%%%%% Model equations %%%%%%%%%%%%%%%%%%%%%%%%%%%% 
\section{The model equations}
\label{sec:modeleq} \indent
 
Let \(\Omega \subset \mathbb{R}^{d}, 2 \le d \le 3\) be a connected open bounded set with a sufficiently smooth boundary \(\partial\Omega\) and denote by \(\mathbf{x}\) the spatial coordinates and by \(t\) the temporal coordinate. The two fluids in \(\Omega\) are considered immiscible and they are contained in the subdomains \(\Omega_{1}(t)\) and \(\Omega_{2}(t)\), respectively, so that \(\overline{\Omega_{1}(t)} \cup \overline{\Omega_{2}(t)} = \overline{\Omega}\). The moving interface between the two fluids is denoted by \(\Gamma(t)\), defined as \(\Gamma(t) = \partial\Omega_{1}(t) \cap \partial\Omega_2(t)\). We consider the classical unsteady, isothermal, incompressible Navier-Stokes equations with gravity, which read as follows \cite{hysing:2009}:
\begin{eqnarray} \label{eq:ns_incomp}
	\rho(\mathbf{x})\left[\frac{\partial \mathbf{u}}{\partial t} + \left(\mathbf{u} \cdot \grad\right)\mathbf{u}\right] &=& -\grad p + \dive\left[2\mu(\mathbf{x}) \mathbf{D}(\mathbf{u})\right] + \rho(\mathbf{x})\mathbf{g} \nonumber \\
	\dive\mathbf{u} &=& 0,   
\end{eqnarray}
for \(\mathbf{x} \in \Omega\), \(t \in (0, T_{f}]\), supplied with suitable initial and boundary conditions, which will be specified in the following. Here \(T_{f}\) is the final time, \(\mathbf{u}\) is the fluid velocity, \(p\) is the pressure, \(\rho\) is the fluid density and \(\mu\) is the dynamic viscosity. We assume that both the density and the viscosity are discontinuous functions
\begin{equation}
	\rho(\mathbf{x}) = 
	\begin{cases}
		\rho_{1} \quad \text{in } \Omega_{1}(t) \\
		\rho_{2} \quad \text{in } \Omega_{2}(t) \\
	\end{cases}
	\qquad \text{and} \qquad
	\mu(\mathbf{x}) = 
	\begin{cases}
		\mu_{1} \quad \text{in } \Omega_{1}(t) \\
		\mu_{2} \quad \text{in } \Omega_{2}(t) \\
	\end{cases}
\end{equation}
with \(\rho_{1}, \rho_{2}, \mu_{1},\) and \(\mu_{2}\) constant values. Moreover, \(\mathbf{g}\) is the gravitational acceleration and \(\mathbf{D}(\mathbf{u})\) denotes the symmetric part of the gradient of the velocity, defined as
\begin{equation}
	\mathbf{D}(\mathbf{u}) = \frac{1}{2}\left[\grad\mathbf{u} + \left(\grad\mathbf{u}\right)^{T}\right].
\end{equation}
In the following, for the sake of simplicity in the notation, we omit the explicit dependence on space and time for the different quantities. Surface tension effects are taken into accounts through the following balance of forces at the interface \(\Gamma\):
\begin{equation}
	\left[\mathbf{u}\right]_{\Gamma} = \mathbf{0} \qquad \left[-p\mathbf{I} + 2\mu \mathbf{D}(\mathbf{u})\right]_{\Gamma}\mathbf{n}_{\Gamma} = \sigma\kappa\mathbf{n}_{\Gamma}, 
\end{equation}
where \(\mathbf{n}_{\Gamma}\) is the outward unit normal to \(\Gamma\), \(\left[\Psi\right]_{\Gamma} = \Psi\rvert_{\Gamma \cap \Omega_{1}} - \Psi\rvert_{\Gamma \cap \Omega_{2}}\) denotes the jump of \(\Psi\) across the interface \(\Gamma\), \(\sigma\) is the constant surface tension coefficient, and \(\kappa = -\dive\mathbf{n}_{\Gamma}\) is the curvature. The first condition implies the continuity of the velocity along \(\Gamma\), whereas the second condition describes the balance of forces at the interface. A common way to handle the term with surface tension is to introduce the following volumetric force \cite{hysing:2009}:
\begin{equation}
	\mathbf{f}_{\sigma} = \sigma\kappa\mathbf{n}_{\Gamma}\delta(\Gamma),
\end{equation} 
where \(\delta(\Gamma)\) is the Dirac delta distribution supported on the interface. Hence, system \eqref{eq:ns_incomp} can be rewritten as follows:
\begin{eqnarray} \label{eq:ns_incomp_surf_tens}
	\rho\left[\frac{\partial \mathbf{u}}{\partial t} + \left(\mathbf{u} \cdot \grad\right)\mathbf{u}\right] &=& -\grad p + \dive\left[2\mu\mathbf{D}(\mathbf{u})\right] + \rho\mathbf{g} + \mathbf{f}_{\sigma} \nonumber \\
	\dive\mathbf{u} &=& 0.   
\end{eqnarray}
A level set approach \cite{osher:2005, sussman:1994} is employed to capture the interface \(\Gamma\). The interface between the two fluids is considered sharp and is described as the zero level set of a smooth function. Hence, the following relation holds:
\begin{equation}\label{eq:levset}
	\frac{\partial\varphi}{\partial t} + \mathbf{u} \cdot \grad\varphi = 0,
\end{equation}
where \(\varphi\) is the level set function. A common choice \cite{sussman:1994} is to consider as level set the signed distance function to \(\Gamma\). In order to fix the notation, we consider \(\varphi < 0\) in \(\Omega_{2}\) and \(\varphi > 0\) in \(\Omega_{1}\). Therefore, we define
\begin{equation}\label{eq:signed_distance}
	\varphi = 
	\begin{cases}
		-\text{dist}(\mathbf{x},\Gamma) \qquad &\text{if } \mathbf{x} \in \Omega_{2} \\
		0 &\text{if } \mathbf{x} \in \Gamma \\
		\text{dist}(\mathbf{x},\Gamma) \qquad &\text{if } \mathbf{x} \in \Omega_{1} \\
	\end{cases}	
\end{equation}
The unit normal vector can be evaluated at each point as follows \cite{dipietro:2006, osher:2005}:
\begin{equation}\label{eq:normal_def}
	\mathbf{n}_{\Gamma} = \frac{\grad\varphi}{\left|\grad\varphi\right|}, \qquad \mathbf{x} \in \Gamma,
\end{equation}
so that \eqref{eq:levset} is equivalent to
\begin{equation}\label{eq:levset_normal}
	\frac{\partial\varphi}{\partial t} + \left(\mathbf{u} \cdot \mathbf{n}_{\Gamma}\right)\left|\grad\varphi\right| = 0.
\end{equation}
Relation \eqref{eq:levset_normal} shows that the deformation of the level set function is due only to the normal component of the velocity. Moreover, we can express the density and the dynamic viscosity through the Heaviside function \(H\)
\begin{eqnarray}
	\rho &=& \rho_{2} + \left(\rho_{1} - \rho_{2}\right)H(\varphi) \\ 
	\mu &=& \mu_{2} + \left(\mu_{1} - \mu_{2}\right)H(\varphi) 
\end{eqnarray}
The whole system of equations reads therefore as follows: 
\begin{eqnarray} \label{eq:ns_incomp_surf_tens_levset}
	\rho\left[\frac{\partial \mathbf{u}}{\partial t} + \left(\mathbf{u} \cdot \grad\right)\mathbf{u}\right] &=& -\grad p + \dive\left[2\mu\mathbf{D}(\mathbf{u})\right] + \rho\mathbf{g} + \mathbf{f}_{\sigma} \nonumber \\
	\dive\mathbf{u} &=& 0 \\
	\frac{\partial\varphi}{\partial t} + \mathbf{u} \cdot \grad\varphi &=& 0. \nonumber  
\end{eqnarray}
System \eqref{eq:ns_incomp_surf_tens_levset} can be rewritten in conservative form. First of all, thanks to the incompressibility constraint \(\dive\mathbf{u} = 0\), we can rewrite \eqref{eq:levset} as
\begin{equation}\label{eq:leveset_cons}
	\frac{\partial\varphi}{\partial t} + \dive\left(\varphi\mathbf{u}\right) = 0.
\end{equation}
Moreover, one can verify that \eqref{eq:levset}, in combination with the incompressibility constraint, implies mass conservation. Indeed, we get
\begin{eqnarray}
	\frac{\partial\rho}{\partial t} + \dive\left(\rho\mathbf{u}\right) = \frac{\partial\rho}{\partial t} + \mathbf{u} \cdot \grad \rho &=& \left(\rho_{1} - \rho_{2}\right)\left(\frac{\partial H(\varphi)}{\partial t} + \mathbf{u} \cdot \grad H(\varphi)\right) \nonumber \\
	&=& \left(\rho_{1} - \rho_{2}\right)\delta(\varphi)\left(\frac{\partial\varphi}{\partial t} + \mathbf{u} \cdot \grad\varphi\right) = 0,
\end{eqnarray}
where we exploited the relation \(\frac{dH(\varphi)}{d\varphi} = \delta(\varphi)\) \cite{salsa:2016}, with \(\delta(\varphi)\) denoting the Dirac delta distribution with support equal to the function \(\varphi\) which implicitly describes the surface. It is appropriate to stress the fact that the differential operators involving the Heaviside function \(H(\varphi)\) have to be intended in a proper distributional sense. Finally, as discussed in \cite{lafaurie:1994}, we can rewrite
\begin{equation}\label{eq:laplace_beltrami}
	\mathbf{f}_{\sigma} = \dive\left[\sigma\left(\mathbf{I} - \mathbf{n}_{\Gamma} \otimes \mathbf{n}_{\Gamma}\right)\delta(\Gamma)\right],
\end{equation}
where, once more, the divergence operator should be intended in a distributional sense. Hence, the conservative form of \eqref{eq:ns_incomp_surf_tens_levset} is
\begin{eqnarray} \label{eq:ns_incomp_surf_tens_levset_cons}
	\frac{\partial\left(\rho\mathbf{u}\right)}{\partial t} + \dive\left(\rho\mathbf{u} \otimes \mathbf{u}\right) &=& -\grad p + \dive\left[2\mu\mathbf{D}(\mathbf{u})\right] + \rho\mathbf{g} + \mathbf{f}_{\sigma} \nonumber \\
	\dive\mathbf{u} &=& 0 \\
	\frac{\partial\varphi}{\partial t} + \dive\left(\varphi\mathbf{u}\right) &=& 0. \nonumber  
\end{eqnarray}
The Continuum Surface Stress (CSS) approach, introduced in \cite{lafaurie:1994}, is employed to treat density, viscosity, and surface tension term. A regularized Heaviside \(H_{\varepsilon}(\varphi)\) is introduced, so as to obtain
\begin{eqnarray}
	\rho \approx \rho_{2} + \left(\rho_{1} - \rho_{2}\right)H_{\varepsilon}(\varphi) \\ 
	\mu \approx \mu_{2} + \left(\mu_{1} - \mu_{2}\right)H_{\varepsilon}(\varphi). 
\end{eqnarray}
It is important at this stage to point out the relation between \(\delta(\Gamma)\) and \(\delta(\varphi)\). As discussed in \cite{estrada:1980}, the following relation holds:
\begin{equation}
	\delta(\Gamma) = \delta(\varphi)\left|\grad\varphi\right|, 
\end{equation} 
so that we can rewrite
\begin{equation}
	\mathbf{f}_{\sigma} = \sigma\kappa\mathbf{n}_{\Gamma}\delta(\varphi)\left|\grad\varphi\right| = \dive\left[\sigma\left(\mathbf{I} - \mathbf{n}_{\Gamma} \otimes \mathbf{n}_{\Gamma}\right)\delta(\varphi)\left|\grad\varphi\right|\right].
\end{equation}
Hence, the CSS approximation of the surface tension term reads as follows:
\begin{equation}
	\mathbf{f}_{\sigma} \approx \dive\left[\sigma\left(\mathbf{I} - \mathbf{n}_{\Gamma} \otimes \mathbf{n}_{\Gamma}\right)\delta_{\varepsilon}(\varphi)\left|\grad\varphi\right|\right] = \dive\left[\sigma\left(\mathbf{I} - \mathbf{n}_{\Gamma} \otimes \mathbf{n}_{\Gamma}\right)\frac{dH_{\varepsilon}}{d\varphi}(\varphi)\left|\grad\varphi\right|\right].
\end{equation}

Since the seminal proposals in \cite{chorin:1968, temam:1969} (see also the review in \cite{guermond:2006}), projection methods have become very popular for the discretization of incompressible Navier-Stokes equations. However, difficulties arise in choosing boundary conditions for the Poisson equation which is to be solved at each time step to compute the pressure. Several strategies have been proposed to deal with this issue and, in particular, high-order pressure boundary conditions were developed in \cite{karniadakis:1991, karniadakis:2005}. This approach, however, requires to rewrite the viscous linear terms as a solenoidal part, which is approximated by an explicit scheme, and as an irrotational part, which is approximate by an implicit scheme of appropriate order \cite{karniadakis:1991}. An alternative that allows to avoid or reduce some of these problems is the so-called artificial compressibility formulation, originally introduced in \cite{chorin:1967} and employed in \cite{bassi:2006, orlando:2022} among many others. In this formulation, the incompressibility constraint is relaxed and a time evolution equation for the pressure is introduced. Hence, we can consider more general boundary conditions for the pressure, and, in particular, full homogenous Neumann boundary conditions can be employed. Moreover, this formulation is consistent with the asymptotic limit of the compressible equations for vanishing Mach number \cite{orlando:2024b}. Hence, the proposed method can be extended to fully-compressible two-phase flows, leading to an asymptotic-preserving scheme, as we aim to show in future work. This kind of approach has been adopted for incompressible flows with variable density, see e.g. \cite{bassi:2022, manzanero:2020}, and we aim here to consider an artificial compressibility formulation for immiscible, isothermal two-phase flows with gravity. The model equations can be therefore rewritten as follows:
\begin{eqnarray} \label{eq:ns_art_comp_surf_tens_levset_cons}
	\frac{\partial\left(\rho\mathbf{u}\right)}{\partial t} + \dive\left(\rho\mathbf{u} \otimes \mathbf{u}\right) &=& -\grad p + \dive\left[2\mu\mathbf{D}(\mathbf{u})\right] + \rho\mathbf{g} + \mathbf{f}_{\sigma} \nonumber \\
	\frac{1}{\rho_{0} c^{2}}\frac{\partial p}{\partial t} + \dive\mathbf{u} &=& 0 \\
	\frac{\partial\varphi}{\partial t} + \dive\left(\varphi\mathbf{u}\right) &=& 0, \nonumber  
\end{eqnarray}
where \(c\) is the artificial speed of sound and \(\rho_{0}\) is a reference density. Finally, since we are relaxing the incompressibility constraint, we consider \eqref{eq:levset} for the level set motion, which is valid for the transport of \(\varphi\) independently of the constraints on the velocity \(\mathbf{u}\). Moreover, this choice is justified by the results reported in \cite{orlando:2023b} for a rising bubble test case, for which a non-conservative formulation leads to less diffusion in the treatment of the interface. Hence, the final form of the system under consideration reads as follows:
\begin{eqnarray} \label{eq:ns_art_comp_surf_tens_levset_fin}
	\frac{\partial\left(\rho\mathbf{u}\right)}{\partial t} + \dive\left(\rho\mathbf{u} \otimes \mathbf{u}\right) &=& -\grad p + \dive\left[2\mu\mathbf{D}(\mathbf{u})\right] + \rho\mathbf{g} + \mathbf{f}_{\sigma} \nonumber \\
	\frac{1}{\rho_{0}c^{2}}\frac{\partial p}{\partial t} + \dive\mathbf{u} &=& 0 \\
	\frac{\partial\varphi}{\partial t} + \mathbf{u} \cdot \grad\varphi &=& 0, \nonumber  
\end{eqnarray}
Before proceeding to describe the time and space discretization schemes, we perform a dimensional analysis to derive a non-dimensional version of system \eqref{eq:ns_art_comp_surf_tens_levset_fin}.

\subsection{Dimensional analysis}
\label{ssec:dim_analysis} 

In this Section, we derive a non-dimensional formulation for system \eqref{eq:ns_art_comp_surf_tens_levset_fin}. We denote with the symbol \(*\) non-dimensional quantities. We introduce a reference length and velocity, denoted by \(L_{ref}\) and \(U_{ref}\), respectively, so as to obtain
\begin{equation}
	\mathbf{x} = L_{ref}\mathbf{x}^{*} \qquad \mathbf{u} = U_{ref}\mathbf{u}^{*} \qquad t = \frac{L_{ref}}{U_{ref}}t^{*}.
\end{equation}
Moreover, we choose as reference density and viscosity those associated to the heavier fluid, which is conventionally considered in \(\Omega_{1}\). For the sake of simplicity, we also assume \(\rho_{0} = \rho_{1}\). The reference pressure \(p_{ref}\) is taken equal to \(p_{ref} = \rho_{1}U_{ref}^{2}\). Hence, we get
\begin{equation}
	\rho = \rho_{1}\rho^{*} \qquad \mu = \mu_{1}\mu^{*} \qquad p = \rho_{1}U_{ref}^{2}p^{*} \qquad \kappa = \frac{1}{L_{ref}}\kappa^{*} \qquad \varphi = L_{ref}\varphi^{*}.
\end{equation}
Introducing the appropriate non-dimensional quantities, we obtain
\begin{eqnarray} \label{eq:ns_art_comp_surf_tens_levset_cons_adim_proof}
	\frac{\rho_{1}U_{ref}^{2}}{L_{ref}} \frac{\partial^{*}\left(\rho^{*}\mathbf{u}^{*}\right)}{\partial^{*}t^{*}} + \frac{\rho_{1}U_{ref}^{2}}{L_{ref}} \grad^{*} \cdot \left(\rho^{*}\mathbf{u}^{*} \otimes \mathbf{u}^{*}\right) &=& -\frac{\rho_{1}U_{ref}^{2}}{L_{ref}} \grad^{*}p^{*} + \frac{\mu_{1}U_{ref}}{L_{ref}} \grad^{*} \cdot \left[2\mu^{*}\mathbf{D}(\mathbf{u}^{*})\right] \nonumber \\
	&-& \rho_{1}\rho^{*}g\mathbf{k} \nonumber \\
	&+& \frac{1}{L_{ref}^{2}} \grad^{*} \cdot \left[\sigma\left(\mathbf{I} - \mathbf{n}_{\Gamma} \otimes \mathbf{n}_{\Gamma}\right)\delta^{*}_{\varepsilon}(\varphi^{*})\left|\grad^{*}\varphi^{*}\right|\right] \nonumber \\
	\frac{\rho_{1}U_{ref}^{3}}{\rho_{1}L_{ref}}
	\frac{1}{c^{2}}\frac{\partial^{*}p^{*}}{\partial^{*}t^{*}} + \frac{U_{ref}}{L_{ref}}\grad^{*} \cdot \mathbf{u}^{*} &=& 0 \\
	U_{ref} \frac{\partial^{*}\varphi^{*}}{\partial^{*}t^{*}} + U_{ref} \mathbf{u}^{*} \cdot \grad^{*}\varphi^{*} &=& 0, \nonumber  
\end{eqnarray}
where \(\mathbf{k}\) is the upward pointing unit vector in the standard Cartesian reference frame. System \eqref{eq:ns_art_comp_surf_tens_levset_cons_adim_proof} reduces to
\begin{eqnarray}
	\frac{\partial^{*}\left(\rho^{*}\mathbf{u}^{*}\right)}{\partial^{*}t^{*}} + \grad^{*} \cdot \left(\rho^{*}\mathbf{u}^{*} \otimes \mathbf{u}^{*}\right) &=& -\grad^{*}p^{*} + \frac{1}{Re} \grad^{*} \cdot \left[2\mu^{*}\mathbf{D}(\mathbf{u}^{*})\right] - \frac{1}{Fr^{2}}\rho^{*}\mathbf{k} \nonumber \\
	&+& \frac{1}{We} \grad^{*} \cdot \left[\left(\mathbf{I} - \mathbf{n}_{\Gamma} \otimes \mathbf{n}_{\Gamma}\right)\delta^{*}_{\varepsilon}(\varphi^{*})\left|\grad^{*}\varphi^{*}\right|\right] \nonumber \\
	M^{2} \frac{\partial^{*}p^{*}}{\partial^{*}t^{*}} + \grad^{*} \cdot \mathbf{u}^{*} &=& 0 \\
	\frac{\partial^{*}\varphi^{*}}{\partial^{*}t^{*}} + \mathbf{u}^{*} \cdot \grad^{*}\varphi^{*} &=& 0, \nonumber  
\end{eqnarray}
where 
\begin{equation}\label{eq:adimensional_parameters}
	Re = \frac{\rho_{1}U_{ref}L_{ref}}{\mu_{1}} \qquad Fr = \frac{U_{ref}}{\sqrt{g L_{ref}}} \qquad We = \frac{\rho_{1}U_{ref}^{2}L_{ref}}{\sigma} \qquad M = \frac{U_{ref}}{c}
\end{equation}
denote the Reynolds number, the Froude number, the Weber number, and the Mach number, respectively. In the following, with a slight abuse of notation, we omit the symbol \(*\) to mark non-dimensional quantities and we consider therefore the following system of equations:
\begin{eqnarray} \label{eq:ns_art_comp_surf_tens_levset_cons_adim}
	\frac{\partial\left(\rho\mathbf{u}\right)}{\partial t} + \dive \left(\rho\mathbf{u} \otimes \mathbf{u}\right) &=& -\grad p + \frac{1}{Re} \dive\left[2\mu \mathbf{D}(\mathbf{u})\right] \nonumber \\
	&-& \frac{1}{Fr^{2}} \rho\mathbf{k} + \frac{1}{We} \dive \left[\left(\mathbf{I} - \mathbf{n}_{\Gamma} \otimes \mathbf{n}_{\Gamma}\right)\delta_{\varepsilon}(\varphi)\left|\grad\varphi\right|\right] \nonumber \\
	M^{2}\frac{\partial p}{\partial t} + \dive \mathbf{u} &=& 0 \\
	\frac{\partial \varphi}{\partial t} + \mathbf{u} \cdot \grad\varphi &=& 0, \nonumber  
\end{eqnarray}
where
\begin{eqnarray}
	\rho &=& \frac{\rho_{2}}{\rho_{1}} + \left(1 - \frac{\rho_{2}}{\rho_{1}}\right)H_{\varepsilon}(\varphi) \\ 
	\mu &=& \frac{\mu_{2}}{\mu_{1}} + \left(1 - \frac{\mu_{2}}{\mu_{1}}\right)H_{\varepsilon}(\varphi). \label{eq:mu_linear}
\end{eqnarray}
	
\subsection{The conservative level set method}
\label{ssec:clsm}

The traditional level set method lacks of volume conservation properties \cite{fedkin:2003}. The conservative level set (CLS) method \cite{olsson:2005, olsson:2007, zahedi:2011} is a popular alternative to add conservation properties to level set schemes. The idea is to replace the signed distance function defined in \eqref{eq:signed_distance} with a regularized Heaviside function:
\begin{equation}\label{eq:conservative_level_set_function}
	\phi(\mathbf{x},t) = \frac{1}{1 + e^{-\varphi(\mathbf{x},t)/\varepsilon}},
\end{equation}
where \(\varepsilon\) helps smoothing the transition of the discontinuous physical properties between the two subdomains and it is also known as interface thickness. Since
\begin{equation}
	\grad\phi = \frac{1}{\varepsilon}\frac{e^{-\varphi/\varepsilon}}{\left(1 + e^{-\varphi/\varepsilon}\right)^2}\grad\varphi
\end{equation}
we can compute the outward unit normal \(\mathbf{n}_\Gamma\) exactly as in \eqref{eq:normal_def}. From definition \eqref{eq:conservative_level_set_function}, it follows that
\begin{equation}
	\Gamma(t) = \left\{\mathbf{x} \in \overline{\Omega} : \phi(\mathbf{x},t) = \frac{1}{2}\right\}.
\end{equation}
This new level set function needs to be reinitialized in order to keep the property of being a regularized Heaviside function \cite{olsson:2007}. This goal is achieved by solving the following PDE \cite{olsson:2005, olsson:2007}:
\begin{equation}\label{eq:reinitialization}
	\frac{\partial\phi}{\partial\tau} + \dive\left(u_{c}\phi\left(1 - \phi\right)\mathbf{n}_{\Gamma}\right) = \dive\left(\beta\varepsilon u_{c}\left(\grad\phi \cdot \mathbf{n}_{\Gamma}\right)\mathbf{n}_{\Gamma}\right), 
\end{equation}
where \(\tau\) is an artificial pseudo-time variable, \(u_{c}\) is an artificial compression velocity, and \(\beta\) is a constant. It is important to notice that \(\mathbf{n}_{\Gamma}\) does not change during the reinizialization procedure, but is computed using the initial value of the level set function. The relation \eqref{eq:reinitialization} has been originally introduced as an intermediate step between the level set advection and the Navier-Stokes equations to keep the shape of the profile \cite{olsson:2005} and to stabilize the advection \cite{olsson:2007}. Two fluxes are considered: a compression flux which acts where \(0 < \phi < 1\) and in normal direction to the interface, represented by \(u_{c}\phi\left(1 - \phi\right)\mathbf{n}_{\Gamma}\), and a diffusion flux, represented by \(\beta\varepsilon u_{c}\left(\grad\phi \cdot \mathbf{n}_{\Gamma}\right)\mathbf{n}_{\Gamma}\). The reinitialization is crucial for the overall stability of the algorithm, but it also introduces errors in the solution \cite{olsson:2007, owkes:2013}. Hence, it is important to avoid unnecessary reinitialization. For this purpose, unlike the formulation proposed e.g. in \cite{olsson:2007} and \cite{owkes:2013}, we introduce the coefficient \(\beta\) to tune the amount of diffusion so as to keep it as small as possible. The choices for the different parameters will be specified in Section \ref{sec:tests}. Finally, we stress the fact that, in this method, we are already using a smooth version of Heaviside function so that
\begin{eqnarray}
	H_\varepsilon &=& \phi \\
	\delta(\Gamma) \approx \frac{dH_\varepsilon}{d\phi}\left|\grad\phi\right| &=& \left|\grad\phi\right|
\end{eqnarray} 

%%%%%%%%%%%%%%%%%%%%%%%%%%%%%%%%%% Time discretization %%%%%%%%%%%%%%%%%%%%%%%%%
\section{The time discretization}
\label{sec:time_disc}

In this Section, we outline the time discretization strategy for system \eqref{eq:ns_art_comp_surf_tens_levset_cons_adim}. Our goal here is to extend the projection method based on the TR-BDF2 scheme developed in \cite{orlando:2022}. We now briefly recall for the convenience of the reader the formulation of the TR-BDF2. This second order implicit method has been originally introduced in \cite{bank:1985} as a combination of the Trapezoidal Rule (or Crank-Nicolson) method
and of the Backward Differentiation Formula method of order 2 (BDF2). Let \(\Delta t = T_{f}/N\) be a discrete time step and \(t^{n} = n\Delta t, n = 0, \dots, N\), be discrete time levels for a generic time dependent problem \(\bm{u}^{\prime} = \mathcal{N}(\bm{u})\). Hence, the incremental form of the TR-BDF2 scheme can be described in terms of two stages, the first one from \(t^{n}\) to \(t^{n+\gamma} = t^{n} + \gamma\Delta t\), and the second one from \(t^{n+\gamma}\) to \(t^{n+1}\), as follows:
\begin{eqnarray}
	\frac{\bm{u}^{n+\gamma} - \bm{u}^{n}}{\gamma\Delta t} &=& \frac{1}{2}\mathcal{N}\left(\bm{u}^{n+\gamma}\right) + \frac{1}{2}\mathcal{N}\left(\bm{u}^{n}\right) \label{eq:trbdf2_stage1} \\
	\frac{\bm{u}^{n+1} - \bm{u}^{n+\gamma}}{\left(1 - \gamma\right)\Delta t} &=& \frac{1}{2 - \gamma}\mathcal{N}\left(\bm{u}^{n+1}\right) + \frac{1 - \gamma}{2\left(2 - \gamma\right)}\mathcal{N}\left(\bm{u}^{n+\gamma}\right) + \frac{1 - \gamma}{2\left(2 - \gamma\right)}\mathcal{N}\left(\bm{u}^{n}\right). \label{eq:trbdf2_stage2}
\end{eqnarray}
Here, \(\bm{u}^{n}\) denotes the approximation at time  \(n = 0, \dots, N\). Notice that, in order to guarantee L-stability, one has to choose \(\gamma = 2 - \sqrt{2}\) \cite{hosea:1996}. We refer to \cite{bonaventura:2017, hosea:1996} for a more exhaustive discussion on the TR-BDF2 method.

We start by considering the equation in system \eqref{eq:ns_art_comp_surf_tens_levset_cons_adim} associated to the level set. In order to avoid a full coupling with the Navier-Stokes equations, we perform a linearization in velocity, so that the first stage for the level set update reads as follows:
\begin{equation}\label{eq:first_stage_levset}
	\frac{\phi^{n+\gamma} - \phi^{n}}{\gamma \Delta t} + \frac{1}{2}\mathbf{u}^{n + \frac{\gamma}{2}} \cdot \grad\phi^{n+\gamma} = -\frac{1}{2}\mathbf{u}^{n + \frac{\gamma}{2}} \cdot \grad\phi^{n},
\end{equation}
where the approximation \(\mathbf{u}^{n + \frac{\gamma}{2}}\) is defined by extrapolation as
\begin{equation}\label{eq:first_stage_extrapolation}
	\mathbf{u}^{n + \frac{\gamma}{2}} = \left(1 + \frac{\gamma}{2\left(1-\gamma\right)}\right)\mathbf{u}^{n} - \frac{\gamma}{2\left(1-\gamma\right)}\mathbf{u}^{n-1}.
\end{equation}	
Following then the projection approach described in \cite{dellarocca:2018, orlando:2022} and applying \eqref{eq:trbdf2_stage1}, the momentum predictor equation for the first stage reads as follows:
\begin{eqnarray} \label{eq:first_stage_momentum}
	&&\frac{\rho^{n+\gamma}\mathbf{u}^{n+\gamma,*} - \rho^{n}\mathbf{u}^{n}}{\gamma\Delta t} \nonumber \\
	&+& \frac{1}{2}\dive\left(\rho^{n+\gamma}\mathbf{u}^{n+\gamma,*} \otimes \mathbf{u}^{n+\frac{\gamma}{2}}\right) - \frac{1}{2}\frac{1}{Re}\dive\left[2\mu^{n+\gamma}\mathbf{D}(\mathbf{u}^{n+\gamma,*})\right] = \nonumber \\
	&-&\frac{1}{2}\dive\left(\rho^{n}\mathbf{u}^{n} \otimes \mathbf{u}^{n+\frac{\gamma}{2}}\right) + \frac{1}{2}\frac{1}{Re}\dive\left[2\mu^{n}\mathbf{D}(\mathbf{u}^{n})\right] - \grad p^n \nonumber \\
	&+& \frac{1}{2}\frac{1}{We}\dive\left[\left(\mathbf{I} - \mathbf{n}_{\Gamma}^{n+\gamma} \otimes \mathbf{n}_{\Gamma}^{n+\gamma}\right)\delta_{\varepsilon}(\phi^{n+\gamma})\left|\grad\phi^{n+\gamma}\right|\right] \\
	&+& \frac{1}{2}\frac{1}{We}\dive\left[\left(\mathbf{I} - \mathbf{n}_{\Gamma}^{n} \otimes\mathbf{n}_{\Gamma}^{n}\right)\delta_{\varepsilon}(\phi^{n})\left|\grad\phi^{n}\right|\right] \nonumber \\
	&-&\frac{1}{2}\frac{1}{Fr^{2}}\rho^{n+\gamma}\mathbf{k} -\frac{1}{2}\frac{1}{Fr^{2}}\rho^{n}\mathbf{k}. \nonumber
\end{eqnarray}
Notice once more that, in order to avoid solving a non-linear system at each time step, \(\mathbf{u}^{n + \frac{\gamma}{2}}\) is employed in the momentum advection terms. We set then \(\delta p^{n+\gamma} = p^{n+\gamma} - p^{n}\) and impose
\begin{eqnarray}\label{eq:first_stage_projection}
	&&\rho^{n+\gamma}\frac{\mathbf{u}^{n+\gamma} - \mathbf{u}^{n+\gamma,*}}{\gamma\Delta t} =-\grad \delta p^{n+\gamma}  \nonumber \\
	&&M^{2}\frac{\delta p^{n+\gamma}}{\gamma\Delta t} + \dive\mathbf{u}^{n+\gamma} = 0.
\end{eqnarray}
Substituting the first equation into the second in \eqref{eq:first_stage_projection}, one obtains the Helmholtz equation
\begin{equation}\label{eq:first_stage_helmholtz}
	M^{2}\frac{\delta p^{n+\gamma}}{\gamma^{2}\Delta t^{2}} - \dive\left(\frac{\grad \delta p^{n+\gamma}}{\rho^{n+\gamma}}\right) = -\frac{1}{\gamma\Delta t} \dive\mathbf{u}^{n+\gamma,*}.
\end{equation}
Once this equation is solved, the final velocity update for the first stage is given by
\begin{equation}
	\mathbf{u}^{n+\gamma} = \mathbf{u}^{n+\gamma,*} - \gamma\Delta t\frac{\grad \delta p^{n+\gamma}}{\rho^{n+\gamma}}.
\end{equation}
The second TR-BDF2 stage is performed in a similar manner applying \eqref{eq:trbdf2_stage2}. We first focus on the level set update:
\begin{equation}\label{eq:second_stage_levset}
	\frac{\phi^{n+1} - \phi^{n+\gamma}}{\left(1-\gamma\right)\Delta t} + a_{33}\mathbf{u}^{n + \frac{3}{2}\gamma} \cdot \grad \phi^{n+1} = - a_{32}\mathbf{u}^{n+\gamma} \cdot \grad\phi^{n+\gamma} - a_{31}\mathbf{u}^{n} \cdot \grad\phi^{n},	
\end{equation}
where
\begin{equation}
	a_{31} = \frac{1-\gamma}{2\left(2-\gamma\right)} \qquad
	a_{32} = \frac{1-\gamma}{2\left(2-\gamma\right)} \qquad
	a_{33} = \frac{1}{2-\gamma}.
\end{equation}
Again, in order to avoid a full coupling with the Navier-Stokes equations, 
an approximation is introduced in the advection term, so that \(\mathbf{u}^{n + \frac{3}{2}\gamma}\) is defined by extrapolation as
\begin{equation}\label{eq:second_stage_extrapolation}
	\mathbf{u}^{n + \frac{3}{2}\gamma} = \left(1 + \frac{1 + \gamma}{\gamma}\right)\mathbf{u}^{n+\gamma} - \frac{1-\gamma}{\gamma}\mathbf{u}^{n}. 
\end{equation}
Then, we define the second momentum predictor:
\begin{eqnarray} \label{eq:second_stage_momentum}
	&&\frac{\rho^{n+1}\mathbf{u}^{n+1,*} - \rho^{n+\gamma}\mathbf{u}^{n+\gamma}}{\left(1-\gamma\right)\Delta t} \nonumber \\
	&+&a_{33}\dive\left(\rho^{n+1}\mathbf{u}^{n+1,*} \otimes \mathbf{u}^{n + \frac{3}{2}\gamma}\right) - a_{33}\frac{1}{Re}\dive\left[2\mu^{n+1}\mathbf{D}(\mathbf{u}^{n+1,*})\right] = \nonumber \\
	&-&a_{32}\dive\left(\rho^{n+\gamma}\mathbf{u}^{n+\gamma} \otimes \mathbf{u}^{n+\gamma}\right) + a_{32}\frac{1}{Re}\dive\left[2\mu^{n+\gamma}\mathbf{D}(\mathbf{u}^{n+\gamma})\right] \nonumber \\
	&-&a_{31}\dive\left(\rho^{n}\mathbf{u}^{n} \otimes \mathbf{u}^{n}\right) + a_{31}\frac{1}{Re}\dive\left[2\mu^{n}\mathbf{D}(\mathbf{u}^{n})\right] - \grad p^{n+\gamma} \nonumber \\
	&+&a_{33}\frac{1}{We}\dive\left[\left(\mathbf{I} - \mathbf{n}_{\Gamma}^{n+1} \otimes \mathbf{n}_{\Gamma}^{n+1}\right)\delta_{\varepsilon}(\phi^{n+1})\left|\grad\phi^{n+1}\right|\right] \\
	&+&a_{32}\frac{1}{We}\dive\left[\left(\mathbf{I} - \mathbf{n}_{\Gamma}^{n+\gamma} \otimes \mathbf{n}_{\Gamma}^{n+\gamma}\right)\delta_{\varepsilon}(\phi^{n+\gamma})\left|\grad\phi^{n+\gamma}\right|\right] \nonumber \\
	&+& a_{31}\frac{1}{We}\dive\left[\left(\mathbf{I} - \mathbf{n}_{\Gamma}^{n} \otimes \mathbf{n}_{\Gamma}^{n}\right)\delta_{\varepsilon}(\phi^{n})\left|\grad\phi^{n}\right|\right] \nonumber \\
	&-&a_{33}\frac{1}{Fr^{2}}\rho^{n+1}\mathbf{k} - a_{32}\frac{1}{Fr^{2}}\rho^{n+\gamma}\mathbf{k} - a_{31}\frac{1}{Fr^{2}}\rho^{n}\mathbf{k}. \nonumber
\end{eqnarray}
Notice that \(\mathbf{u}^{n + \frac{3}{2}\gamma}\) is employed in the non-linear momentum advection term. We set then $\delta p^{n+1} =p^{n+1} - p^{n+\gamma} $ and impose
\begin{eqnarray}\label{eq:second_stage_projection}
	&&\rho^{n+1}\frac{\mathbf{u}^{n+1} - \mathbf{u}^{n+1,*}}{(1-\gamma)\Delta t} =-\grad \delta p^{n+1} \nonumber \\
	&&M^{2}\frac{\delta p^{n+1}}{(1-\gamma)\Delta t} + \dive\mathbf{u}^{n+1} = 0.
\end{eqnarray}
Substituting the first equation into the second in \eqref{eq:second_stage_projection}, one obtains the Helmholtz equation
\begin{equation}\label{eq:second_stage_helmholtz}
	M^{2}\frac{\delta p^{n+1}}{(1-\gamma)^{2}\Delta t^{2}} -\dive\left(\frac{\grad \delta p^{n+1}}{\rho^{n+1}}\right) = -\frac{1}{(1-\gamma)\Delta t} \dive\mathbf{u}^{n+1,*}.
\end{equation}
The final velocity update then reads as follows:
\begin{equation}
	\mathbf{u}^{n+1} = \mathbf{u}^{n+1,*} - (1 - \gamma)\Delta t\frac{\grad \delta p^{n+1}}{\rho^{n+1}}.
\end{equation}
Finally, we focus on the reinitialization procedure described in Equation \ref{eq:reinitialization}, which is performed after each stage of the level set update and before computing the momentum predictor. We consider an implicit treatment of the diffusion term \(\dive\left(\beta\varepsilon u_{c}\left(\grad\phi\cdot\mathbf{n}_\Gamma\right)\mathbf{n}_\Gamma\right)\) and a semi-implicit treatment of the compression term \(u_{c}\dive\left(\phi\left(1 - \phi\right)\mathbf{n}_{\Gamma}\right)\). Hence, the semi-discrete formulation reads as follows:
\begin{equation}\label{eq:reinitialization_time}
	\frac{\phi^{k+1,*} - \phi^{k,*}}{\Delta \tau} + \dive\left(u_{c}\phi^{k+1,*}\left(1 - \phi^{k,*}\right)\mathbf{n}_{\Gamma}\right) = \dive\left(\beta\varepsilon u_{c}\left(\grad\phi^{k+1,*} \cdot \mathbf{n}_{\Gamma}\right)\mathbf{n}_{\Gamma}\right), 
\end{equation}
where \(\Delta\tau\) is the pseudo time step. Moreover, \(\phi^{0,*} = \phi^{n+\gamma}\) after the first TR-BDF2 stage and \(\phi^{0,*} = \phi^{n+1}\) after the second TR-BDF2 stage. We recall once more that \(\mathbf{n}_{\Gamma} = \frac{\grad\phi^{0,*}}{\left|\grad\phi^{0,*}\right|}\) and it does not change during the reinitialization. Following \cite{owkes:2013}, we define the total reinitialization time \(\tau_{fin}\) as a fraction of the time step \(\Delta t\), namely
\begin{equation}
	\tau_{fin} = \eta \Delta t.
\end{equation}
\(\eta = 0\) corresponds to no reinitialization, whereas \(\eta = 1\) yields an amount of reinitialization which can modify the values of level set function of the same order of magnitude of which they have been modified during the previous advection step. For most applications, \(\eta \approx 0.5\) seems to provide an appropriate amount of reinitialization \cite{owkes:2013}. A pseudo time step such that two to five reinitialization steps are performed typically ensures stable solutions and leads to the updated level set function \cite{kronbichler:2018}.

%%%%%%%%%%%%%%%%%%%%%%%%%%%% Spatial discretization %%%%%%%%%%%%%%%%%%%%%%%%%%%%%
\section{The spatial discretization}
\label{sec:space_disc} \indent
 
For the spatial discretization, we consider discontinuous finite element approximations. We consider a decomposition of the domain \(\Omega\) into a family of hexahedra \(\mathcal{T}_{h}\) (quadrilaterals in the two-dimensional case) and denote each element by \(K\). The skeleton \(\mathcal{E}\) denotes the set of all element faces and \(\mathcal{E} = \mathcal{E}^{I} \cup \mathcal{E}^{B}\), where \(\mathcal{E}^{I}\) is the subset of interior faces and \(\mathcal{E}^{B}\) is the subset of boundary faces. Suitable jump and average operators can then be defined as customary for finite element discretizations. A face \(e \in \mathcal{E}^{I}\) shares two elements that we denote by \(K^{+}\) with outward unit normal \(\mathbf{n}^{+}\) and \(K^{-}\) with outward unit normal \(\mathbf{n}^{-}\), whereas for a face \(e \in \mathcal{E}^{B}\) we denote by \(\mathbf{n}\) the outward unit normal. For a scalar function \(\Psi\) the jump is defined as
\begin{equation}
	\left[\left[\Psi\right]\right] = \Psi^{+}\mathbf{n}^{+} + \Psi^{-}\mathbf{n}^{-} \quad \text{if } e \in \mathcal{E}^{I} \qquad \left[\left[\Psi\right]\right] = \Psi\mathbf{n} \quad \text{if } e \in \mathcal{E}^{B}.
\end{equation}
The average is defined as
\begin{equation}
	\left\{\left\{\Psi\right\}\right\} = \frac{1}{2}\left(\Psi^{+} + \Psi^{-}\right) \quad \text{if } e \in \mathcal{E}^{I} \qquad \left\{\left\{\Psi\right\}\right\} = \Psi \quad \text{if } e \in \mathcal{E}^{B}.	
\end{equation}
Similar definitions apply for a vector function \(\boldsymbol{\Psi}\):
\begin{align}
	&\left[\left[\boldsymbol{\Psi}\right]\right] = \boldsymbol{\Psi}^{+} \cdot \mathbf{n}^{+} + \boldsymbol{\Psi}^{-} \cdot \mathbf{n}^{-} \quad \text{if } e \in \mathcal{E}^{I} \qquad 
	\left[\left[\boldsymbol{\Psi}\right]\right] = \boldsymbol{\Psi} \cdot \mathbf{n} \quad \text{if } e \in \mathcal{E}^{B} \\
	&\left\{\left\{\boldsymbol{\Psi}\right\}\right\} = \frac{1}{2}\left(\boldsymbol{\Psi}^{+} + \boldsymbol{\Psi}^{-}\right) \quad \text{if } e \in \mathcal{E}^{I} \qquad \left\{\left\{\boldsymbol{\Psi}\right\}\right\} = \boldsymbol{\Psi} \quad \text{if } e \in \mathcal{E}^{B}.
\end{align}
For vector functions, it is also useful to define a tensor jump as:
\begin{equation}
	\left<\left<\boldsymbol{\Psi}\right>\right> = \boldsymbol{\Psi}^{+} \otimes \mathbf{n}^{+} + \boldsymbol{\Psi}^{-} \otimes \mathbf{n}^{-} \quad \text{if }\Gamma \in \mathcal{E}^{I} 
	\qquad \left<\left<\boldsymbol{\Psi}\right>\right> = \boldsymbol{\Psi} \otimes \mathbf{n} \quad \text{if }\Gamma \in \mathcal{E}^{B}.	
\end{equation}
We now introduce the following finite element spaces:
\[Q_{k} = \left\{v \in L^2(\Omega) : v\rvert_{K} \in \mathbb{Q}_{k} \quad \forall K \in \mathcal{T}_{h}\right\}\] 
and
\[\mathbf{Q}_{k} = \left[Q_{k}\right]^d,\]
where \(\mathbb{Q}_{k}\) is the space of polynomials of degree \(k\) in each coordinate direction. Considering the well-posedness analyses in \cite{schotzau:2003, toselli:2002}, the finite element spaces that will be used for the discretization of velocity and pressure are \(\mathbf{V}_{h} = \mathbf{Q}_{k}\) and \(W_{h} = Q_{k-1} \cap L^2_{0}(\Omega)\), respectively, where \(k \ge 2\). For what concerns the level set function, we consider instead \(X_{h} = Q_{r}\) with \(r \ge 2\), so that its gradient is at least a piecewise linear polynomial. We then denote by \(\psi_{i}(\mathbf{x})\) the basis functions for the finite element spaces associated to the scalar variable, i.e. \(W_{h}\) and \(X_{h}\), and by \(\boldsymbol{\psi}_{i}(\mathbf{x})\) the basis functions for the space \(V_{h}\), the finite element space chosen for the discretization of the velocity. Hence, we get
\begin{equation}\label{eq:fe_representations}
	\mathbf{u} \approx \sum_{j=1}^{\text{dim}(\mathbf{V}_{h})} u_{j}(t)\boldsymbol{\psi}_{j}(\mathbf{x}) \qquad p \approx \sum_{j=1}^{\text{dim}(W_{h})} p_{j}(t)\psi_{j}(\mathbf{x}) \qquad \phi \approx \sum_{j=1}^{\text{dim}(X_{h})}\phi_{j}(t)\psi_{j}(\mathbf{x})
\end{equation} 
The shape functions correspond to the products of Lagrange interpolation polynomials for the support points of \((k+1)\)-order Gauss-Lobatto quadrature rule in each coordinate direction. Finally, as custom for discontinuous finite elements, we employ a weak imposition of the boundary conditions \cite{arnold:2002}. More specifically, we deduce suitable exterior values from the boundary conditions, employing the so-called mirror principle. For Dirichlet boundary conditions, we set (e.g. for the velocity) \(\mathbf{u}^{-} = - \mathbf{u}^{+} + 2\mathbf{u}_{D}\) with \(\left[\grad\mathbf{u}^{+}\right]\mathbf{n} = \left[\grad\mathbf{u}^{-}\right]\mathbf{n}\) and \(\mathbf{u}_{D}\) denoting the Dirichlet value. On the other hand, for Neumann boundary conditions, we set (e.g. for the pressure) \(\grad p^{-} \cdot \mathbf{n} = -\grad p^{+} \cdot \mathbf{n} + 2p_{N}\) with \(p^{+} = p^{-}\) and \(p_{N}\) denoting the Neumann value. Given these definitions, the weak formulation of the level set update for the first stage is obtained multiplying equation \eqref{eq:first_stage_levset} by a test function \(w \in X_{h}\):
\begin{eqnarray}
	&&\sum_{K \in \mathcal{T}_{h}}\int_{K} \frac{\phi^{n+\gamma}}{\gamma \Delta t}w d\Omega + \frac{1}{2}\sum_{K \in \mathcal{T}_{h}}\int_{K} \mathbf{u}^{n+\frac{\gamma}{2}} \cdot \grad\phi^{n+\gamma} w d\Omega \nonumber \\
	&+& \frac{1}{2}\sum_{e \in \mathcal{E}} \int_{e} \left\{\left\{\phi^{n+\gamma}\mathbf{u}^{n + \frac{\gamma}{2}}\right\}\right\} \cdot \left[\left[w\right]\right]d\Sigma - \frac{1}{2}\sum_{e \in \mathcal{E}} \int_{e} \left\{\left\{\mathbf{u}^{n + \frac{\gamma}{2}}\right\}\right\} \cdot \left[\left[\phi^{n + \gamma}w\right]\right]d\Sigma \nonumber \\
	&+& \frac{1}{2}\sum_{e \in \mathcal{E}} \int_{e} \frac{\lambda^{n + \frac{\gamma}{2}}}{2}\left[\left[\phi^{n + \gamma}\right]\right] \cdot \left[\left[w\right]\right]d\Sigma \nonumber \\
	&=&\sum_{K \in \mathcal{T}_{h}}\int_{K} \frac{\phi^{n}}{\gamma \Delta t}w d\Omega - \frac{1}{2}\sum_{K \in \mathcal{T}_{h}}\int_{K} \mathbf{u}^{n + \frac{\gamma}{2}} \cdot \grad\phi^{n} w d\Omega \\
	&-& \frac{1}{2}\sum_{e \in \mathcal{E}} \int_{e} \left\{\left\{\phi^{n}\mathbf{u}^{n + \frac{\gamma}{2}}\right\}\right\} \cdot \left[\left[w\right]\right]d\Sigma - \frac{1}{2}\sum_{e \in \mathcal{E}} \int_{e} \left\{\left\{\mathbf{u}^{n + \frac{\gamma}{2}}\right\}\right\} \cdot \left[\left[\phi^{n}w\right]\right]d\Sigma \nonumber \\
	&-& \frac{1}{2}\sum_{e \in \mathcal{E}} \int_{e} \frac{\lambda^{n + \frac{\gamma}{2}}}{2}\left[\left[\phi^{n}\right]\right] \cdot \left[\left[w\right]\right]d\Sigma, \nonumber
\end{eqnarray}
where
\begin{equation}
	\lambda^{n + \frac{\gamma}{2}} = \max\left(\left|\left(\mathbf{u}^{n + \frac{\gamma}{2}}\right)^{+} \cdot \mathbf{n}^{+}\right|, \left|\left(\mathbf{u}^{n + \frac{\gamma}{2}}\right)^{-} \cdot \mathbf{n}^{-}\right|\right).
\end{equation}
Following \cite{bassi:1997}, the numerical approximation of the non-conservative term is based on a double integration by parts. The algebraic form can be obtained taking \(w = \psi_{i}, i = 1,\dots,\text{dim}(X_{h})\) and exploiting the representation in \eqref{eq:fe_representations}, so as to obtain in compact form
\begin{equation}
	\left(\frac{1}{\gamma \Delta t}\mathbf{M}_{\phi} + \frac{1}{2}\mathbf{A}_{\phi}^{n+\gamma}\right)\makebox{\large\ensuremath{\boldsymbol{\phi}}}^{n+\gamma} = \mathbf{F}_{\phi}^{n},
\end{equation}
where \(\makebox{\large\ensuremath{\boldsymbol{\phi}}}^{n+\gamma}\) denotes the vector of the degrees of freedom associated to the level set. Moreover, we have set
\begin{eqnarray}
	\mathbf{M}_{\phi_{ij}} &=& \sum_{K \in \mathcal{T}_{h}}\int_{K} \psi_{j}\psi_{i} d\Omega \\
	\mathbf{A}_{\phi_{ij}}^{n+\gamma} &=& \sum_{K \in \mathcal{T}_{h}}\int_{K} \mathbf{u}^{n+\frac{\gamma}{2}} \cdot \grad \psi_{j}\psi_{i} d\Omega \nonumber \\
	&+& \sum_{e \in \mathcal{E}} \int_{e} \left\{\left\{\mathbf{u}^{n + \frac{\gamma}{2}}\psi_{j}\right\}\right\} \cdot \left[\left[\psi_{i}\right]\right]d\Sigma - \sum_{e \in \mathcal{E}} \int_{e} \left\{\left\{\mathbf{u}^{n + \frac{\gamma}{2}}\right\}\right\} \cdot \left[\left[\psi_{j}\psi_{i}\right]\right]d\Sigma \nonumber \\
	&+& \sum_{e \in \mathcal{E}} \int_{e} \frac{\lambda^{n + \frac{\gamma}{2}}}{2} \left[\left[\psi_{j}\right]\right] \cdot \left[\left[\psi_{i}\right]\right]d\Sigma 
\end{eqnarray}
and 
\begin{eqnarray}
	\mathbf{F}_{\phi}^{n} &=& \sum_{K \in \mathcal{T}_{h}}\int_{K} \frac{\phi^{n}}{\gamma \Delta t}\psi_{i} d\Omega + \frac{1}{2}\sum_{K \in \mathcal{T}_{h}}\int_{K} \mathbf{u}^{n + \frac{\gamma}{2}} \cdot \grad\phi^{n} \psi_{i} d\Omega \\
	&-& \frac{1}{2}\sum_{e \in \mathcal{E}} \int_{e} \left\{\left\{\phi^{n}\mathbf{u}^{n + \frac{\gamma}{2}}\right\}\right\} \cdot \left[\left[\psi_{i}\right]\right]d\Sigma + \frac{1}{2} \sum_{e \in \mathcal{E}} \int_{e} \left\{\left\{\mathbf{u}^{n + \frac{\gamma}{2}}\right\}\right\} \cdot \left[\left[\phi^{n}\psi_{i}\right]\right]d\Sigma \nonumber \\
	&-& \frac{1}{2}\sum_{e \in \mathcal{E}} \int_{e} \frac{\lambda^{n + \frac{\gamma}{2}}}{2}\left[\left[\phi^{n}\right]\right] \cdot \left[\left[\psi_{i}\right]\right]d\Sigma. \nonumber
\end{eqnarray}
Consider now the variational formulation for equation \eqref{eq:first_stage_momentum}. Take \(\mathbf{v} \in \mathbf{V}_{h}\) so as to obtain after integration by parts
\begin{eqnarray}
	&&\sum_{K \in \mathcal{T}_{h}} \int_{K} \frac{1}{\gamma \Delta t} \rho^{n+\gamma}\mathbf{u}^{n+\gamma,*} \cdot \mathbf{v} d\Omega - \frac{1}{2} \sum_{K \in \mathcal{T}_{h}} \int_{K} \rho^{n+\gamma}\mathbf{u}^{n+\gamma,*} \otimes \mathbf{u}^{n + \frac{\gamma}{2}} : \grad\mathbf{v} d\Omega \nonumber \\
	&+&\frac{1}{2} \sum_{e \in \mathcal{E}} \int_{e} \left\{\left\{\rho^{n+\gamma}\mathbf{u}^{n+\gamma,*} \otimes \mathbf{u}^{n + \frac{\gamma}{2}}\right\}\right\} : \left<\left<\mathbf{v}\right>\right> d\Sigma + \frac{1}{2} \sum_{e \in \mathcal{E}} \int_{e} \frac{\lambda^{n + \frac{\gamma}{2}}}{2}\left<\left<\rho^{n+\gamma}\mathbf{u}^{n + \frac{\gamma}{2}}\right>\right> : \left<\left<\mathbf{v}\right>\right> d\Sigma \nonumber \\
	&+&\frac{1}{2Re} \sum_{K \in \mathcal{T}_{h}} \int_{K} 2\mu^{n+\gamma}\mathbf{D}(\mathbf{u}^{n+\gamma,*}) : \grad\mathbf{v}d\Omega - \frac{1}{2Re} \sum_{e \in \mathcal{E}} \int_{e} \left\{\left\{2\mu^{n+\gamma}\mathbf{D}(\mathbf{u}^{n+\gamma,*})\right\}\right\} : \left<\left<\mathbf{v}\right>\right> d\Sigma \nonumber \\
	&-&\frac{1}{2Re} \sum_{e \in \mathcal{E}} \int_{e} \left<\left<\mathbf{u}^{n+\gamma,*}\right>\right> : \left\{\left\{2\mu^{n+\gamma}\mathbf{D}(\mathbf{v})\right\}\right\} d\Sigma \nonumber \\
	&+&\frac{1}{2Re} \sum_{e \in \mathcal{E}} \int_{e} C_{u} \left\{\left\{\mu^{n+\gamma}\right\}\right\}_{H}\left<\left<\mathbf{u}^{n+\gamma,*}\right>\right> : \left<\left<\mathbf{v}\right>\right>d\Sigma \nonumber \\
	&=&\sum_{K \in \mathcal{T}_{h}} \int_{K} \frac{1}{\gamma \Delta t} \rho^{n}\mathbf{u}^{n} \cdot \mathbf{v} d\Omega + \frac{1}{2} \sum_{K \in \mathcal{T}_{h}} \int_{K} \rho^{n}\mathbf{u}^{n} \otimes \mathbf{u}^{n + \frac{\gamma}{2}} : \grad\mathbf{v} d\Omega \\
	&-&\frac{1}{2} \sum_{e \in \mathcal{E}} \int_{e} \left\{\left\{\rho^{n}\mathbf{u}^{n} \otimes \mathbf{u}^{n + \frac{\gamma}{2}}\right\}\right\} : \left<\left<\mathbf{v}\right>\right> d\Sigma - \frac{1}{2} \int_{e} \frac{\lambda^{n+\frac{\gamma}{2}}}{2}\left<\left<\rho^{n}\mathbf{u}^{n}\right>\right> : \left<\left<\mathbf{v}\right>\right>d\Sigma \nonumber \\
	&-&\frac{1}{2Re} \sum_{K \in \mathcal{T}_{h}} \int_{K} 2\mu^{n}\mathbf{D}(\mathbf{u}^{n}) : \grad\mathbf{v}d\Omega + \frac{1}{2Re} \sum_{e \in \mathcal{E}} \int_{e} \left\{\left\{2\mu^{n}\mathbf{D}(\mathbf{u}^{n})\right\}\right\} : \left<\left<\mathbf{v}\right>\right> d\Sigma \nonumber \\
	&+&\sum_{K \in \mathcal{T}_h} \int_{K} p^{n} \dive \mathbf{v} d\Omega - \sum_{e \in \mathcal{E}} \int_{e} \left\{\left\{p^{n}\right\}\right\}\left[\left[\mathbf{v}\right]\right] d\Sigma \nonumber \\
	&-&\frac{1}{2 Fr^{2}} \sum_{K \in \mathcal{T}_{h}} \int_{K} \rho^{n+\gamma}\mathbf{k} \cdot \mathbf{v}d\Omega - \frac{1}{2 Fr^{2}} \sum_{K \in \mathcal{T}_{h}} \int_{K} \rho^{n}\mathbf{k} \cdot \mathbf{v}d\Omega \nonumber \\
	&-&\frac{1}{2 We} \sum_{K \in \mathcal{T}_{h}} \int_{K} \left(\mathbf{I} - \mathbf{n}_{\Gamma}^{n+\gamma} \otimes \mathbf{n}_{\Gamma}^{n+\gamma}\right)\delta_{\varepsilon}(\phi^{n+\gamma})\left|\grad\phi^{n+\gamma}\right| : \grad\mathbf{v} d\Omega \nonumber \\ 
	&+&\frac{1}{2 We} \sum_{e \in \mathcal{E}} \int_{e} \left\{\left\{\left(\mathbf{I} - \mathbf{n}_{\Gamma}^{n+\gamma} \otimes \mathbf{n}_{\Gamma}^{n+\gamma}\right)\delta_{\varepsilon}(\phi^{n+\gamma})\left|\grad\phi^{n+\gamma}\right|\right\}\right\} : \left<\left<\mathbf{v}\right>\right> d\Sigma \nonumber \\
	&-&\frac{1}{2 We} \sum_{K \in \mathcal{T}_{h}} \int_{K} \left(\mathbf{I} - \mathbf{n}_{\Gamma}^{n} \otimes \mathbf{n}_{\Gamma}^{n}\right)\delta_{\varepsilon}(\phi^{n})\left|\grad\phi^{n}\right| : \grad\mathbf{v} d\Omega \nonumber \\ 
	&+&\frac{1}{2 We} \sum_{e \in \mathcal{E}} \int_{e} \left\{\left\{\left(\mathbf{I} - \mathbf{n}_{\Gamma}^{n} \otimes \mathbf{n}_{\Gamma}^{n}\right)\delta_{\varepsilon}(\phi^{n})\left|\grad\phi^{n}\right|\right\}\right\} : \left<\left<\mathbf{v}\right>\right> d\Sigma, \nonumber
\end{eqnarray}
where
\begin{equation}
	\left\{\left\{\mu^{n+\gamma}\right\}\right\}_{H} = \frac{2}{\frac{1}{\mu^{n+\gamma,+}} + \frac{1}{\mu^{n+\gamma,-}}}.
\end{equation}
Here, following e.g. \cite{antonietti:2012}, we employ the harmonic average of the viscosity coefficient for the penalization term. Notice that the approximation of the advection term employs an upwind flux, whereas the approximation of the diffusion term is based on the Symmetric Interior Penalty (SIP) \cite{arnold:1982}. Notice also that no penalization terms have been introduced for the variables computed at previous time steps in the diffusion terms. Following \cite{fehn:2019, orlando:2022}, we set for each face \(e\) of a cell \(K\) 
\begin{equation}
	\sigma^{\mathbf{u}}_{e,K} = \left(k + 1 \right)^2\frac{\text{diam}(e)}{\text{diam}(K)}
\end{equation}
and we define the penalization constant for the SIP method as 
\begin{equation}
	C_{u} = \frac{1}{2}\left(\sigma^{\mathbf{u}}_{e,K^{+}} + \sigma^{\mathbf{u}}_{e,K^{-}}\right) \quad \text{if } e \in \mathcal{E}^{I}, \qquad C_{u} = \sigma^{\mathbf{u}}_{e,K} \quad \text{if } e \in \mathcal{E}^{B}. 
\end{equation}
Finally, we stress the fact that a centered flux has been employed for the surface tension terms. The algebraic formulation is then computed considering \(\mathbf{v} = \boldsymbol{\psi}_{i}, i=1, \dots, \text{dim}(\mathbf{V}_{h})\) and the representation in \eqref{eq:fe_representations} for the velocity. Hence, we obtain
\begin{equation}
	\left(\frac{1}{\gamma\Delta t}\mathbf{M}_{\mathbf{u}}^{n+\gamma} + \frac{1}{2Re}\mathbf{A}_{\mathbf{u}}^{n+\gamma} + \frac{1}{2}\mathbf{C}_{\mathbf{u}}^{n+\gamma}\right)\mathbf{U}^{n+\gamma,*} = \mathbf{F}_{u}^{n},
\end{equation}
where \(\mathbf{U}^{n+\gamma,*}\) denotes the vector of degrees of freedom for the velocity. Moreover, we have set
\begin{eqnarray}
	\mathbf{M}_{\mathbf{u}_{ij}}^{n+\gamma} &=& \sum_{K \in \mathcal{T}_{h}} \int_{K} \rho^{n+\gamma}\boldsymbol{\psi}_{j}\boldsymbol{\psi}_{i} d\Omega \\
	\mathbf{C}_{\mathbf{u}_{ij}}^{n+\gamma} &=& -\sum_{K \in \mathcal{T}_{h}} \int_{K} \rho^{n+\gamma}\boldsymbol{\psi}_{j} \otimes \mathbf{u}^{n + \frac{\gamma}{2}} : \grad\boldsymbol{\psi}_{i} d\Omega \nonumber \\
	&+& \sum_{e \in \mathcal{E}} \int_{e} \left\{\left\{\rho^{n+\gamma}\boldsymbol{\psi}_{j} \otimes \mathbf{u}^{n + \frac{\gamma}{2}}\right\}\right\} : \left<\left<\boldsymbol{\psi}_{j}\right>\right> d\Sigma	\\
	&+& \sum_{e \in \mathcal{E}} \int_{e} \frac{\lambda^{n + \frac{\gamma}{2}}}{2}\left<\left<\rho^{n+\gamma}\boldsymbol{\psi}_{j}\right>\right> : \left<\left<\boldsymbol{\psi}_{i}\right>\right> d\Sigma \nonumber \\
	\mathbf{A}_{\mathbf{u}_{ij}}^{n+\gamma} &=& \sum_{K \in \mathcal{T}_{h}} \int_{K} 2\mu^{n+\gamma}\mathbf{D}\left(\boldsymbol{\psi}_{j}\right) : \grad\boldsymbol{\psi}_{i} d\Omega \nonumber \\
	&-& \sum_{e \in \mathcal{E}} \int_{e} \left\{\left\{2\mu^{n+\gamma}\mathbf{D}(\boldsymbol{\psi}_{j})\right\}\right\} : \left<\left<\boldsymbol{\psi}_{i}\right>\right>d\Sigma \nonumber \\
	&-& \sum_{e \in \mathcal{E}} \int_{e} \left<\left<\boldsymbol{\psi}_{j}^{n+\gamma,*}\right>\right> : \left\{\left\{2\mu^{n+\gamma}\mathbf{D}(\boldsymbol{\psi}_{i})\right\}\right\}d\Sigma \nonumber \\
	&+& \sum_{e \in \mathcal{E}} \int_{e} C_{u}\left\{\left\{\mu^{n+\gamma}\right\}\right\}_{H} \left<\left<\boldsymbol{\psi}_{j}\right>\right> : \left<\left<\boldsymbol{\psi}_{i}\right>\right>d\Sigma
\end{eqnarray}
and
\begin{eqnarray}
	\mathbf{F}_{\mathbf{u}}^{n} &=& \sum_{K \in \mathcal{T}_{h}} \int_{K} \frac{1}{\gamma \Delta t} \rho^{n}\mathbf{u}^{n} \cdot \boldsymbol{\psi}_{i} d\Omega + \frac{1}{2} \sum_{K \in \mathcal{T}_{h}} \int_{K} \rho^{n}\mathbf{u}^{n} \otimes \mathbf{u}^{n + \frac{\gamma}{2}} : \grad\boldsymbol{\psi}_{i} d\Omega \nonumber \\
	&-& \frac{1}{2} \sum_{e \in \mathcal{E}} \int_{e} \left\{\left\{\rho^{n}\mathbf{u}^{n} \otimes \mathbf{u}^{n + \frac{\gamma}{2}}\right\}\right\} : \left<\left<\boldsymbol{\psi}_{i}\right>\right> d\Sigma \nonumber \\
	&-& \frac{1}{2} \sum_{e \in \mathcal{E}} \int_{e} \frac{\lambda^{n+\frac{\gamma}{2}}}{2}\left<\left<\rho^{n}\mathbf{u}^{n}\right>\right> : \left<\left<\boldsymbol{\psi}_{i}\right>\right> d\Sigma \nonumber \\
	&-& \frac{1}{2Re} \sum_{K \in \mathcal{T}_{h}} \int_{K} 2\mu^{n}\mathbf{D}(\mathbf{u}^{n}) : \grad\boldsymbol{\psi}_{i}d\Omega + \frac{1}{2Re} \sum_{e \in \mathcal{E}} \int_{e} \left\{\left\{2\mu^{n}\mathbf{D}(\mathbf{u}^{n})\right\}\right\} : \left<\left<\boldsymbol{\psi}_{i}\right>\right> d\Sigma \nonumber \\
	&+& \sum_{K \in \mathcal{T}_h} \int_{K} p^{n} \dive \boldsymbol{\psi}_{i} d\Omega - \sum_{e \in \mathcal{E}} \int_{e} \left\{\left\{p^{n}\right\}\right\}\left[\left[\boldsymbol{\psi}_{i}\right]\right] d\Sigma \\
	&-& \frac{1}{2 Fr^{2}} \sum_{K \in \mathcal{T}_{h}} \int_{K} \rho^{n+\gamma}\mathbf{k} \cdot \boldsymbol{\psi}_{i}d\Omega - \frac{1}{2 Fr^{2}} \sum_{K \in \mathcal{T}_{h}} \int_{K} \rho^{n}\mathbf{k} \cdot \boldsymbol{\psi}_{i}d\Omega \nonumber \\
	&-& \frac{1}{2 We} \sum_{K \in \mathcal{T}_{h}} \int_{K} \left(\mathbf{I} - \mathbf{n}_{\Gamma}^{n+\gamma} \otimes \mathbf{n}_{\Gamma}^{n+\gamma}\right)\delta_{\varepsilon}(\phi^{n+\gamma})\left|\grad\phi^{n+\gamma}\right| : \grad\boldsymbol{\psi}_{i}d\Omega \nonumber \\ 
	&+& \frac{1}{2 We} \sum_{e \in \mathcal{E}} \int_{e} \left\{\left\{\left(\mathbf{I} - \mathbf{n}_{\Gamma}^{n+\gamma} \otimes \mathbf{n}_{\Gamma}^{n+\gamma}\right)\delta_{\varepsilon}(\phi^{n+\gamma})\left|\grad\phi^{n+\gamma}\right|\right\}\right\} : \left<\left<\boldsymbol{\psi}_{i}\right>\right> d\Sigma \nonumber \\
	&-& \frac{1}{2 We} \sum_{K \in \mathcal{T}_{h}} \int_{K} \left(\mathbf{I} - \mathbf{n}_{\Gamma}^{n} \otimes \mathbf{n}_{\Gamma}^{n}\right)\delta_{\varepsilon}(\phi^{n})\left|\grad\phi^{n}\right| : \grad\boldsymbol{\psi}_{i} d\Omega \nonumber \\ 
	&+& \frac{1}{2 We} \sum_{e \in \mathcal{E}} \int_{e} \left\{\left\{\left(\mathbf{I} - \mathbf{n}_{\Gamma}^{n} \otimes \mathbf{n}_{\Gamma}^{n}\right)\delta_{\varepsilon}(\phi^{n})\left|\grad\phi^{n}\right|\right\}\right\} : \left<\left<\boldsymbol{\psi}_{i}\right>\right> d\Sigma. \nonumber	
\end{eqnarray}
For what concerns the projection step, we apply again the SIP method. We multiply \eqref{eq:first_stage_helmholtz} by a test function \(q \in Q_{h}\), we apply Green's theorem and we get
\begin{eqnarray}
	&&\sum_{K \in \mathcal{T}_{h}} \int_{K} \frac{M^{2}}{\gamma^{2}\Delta t^{2}} \delta p^{n+\gamma} q d\Omega + \sum_{K \in \mathcal{T}_{h} }\int_{K} \frac{\grad \delta p^{n+\gamma}}{\rho^{n+\gamma}} \cdot \grad q d\Omega \nonumber \\
	&-& \sum_{e \in \mathcal{E}} \int_{e} \left\{\left\{\frac{\grad \delta p^{n+\gamma}}{\rho^{n+\gamma}}\right\}\right\} \cdot \left[\left[q\right]\right]d\Sigma - \sum_{e \in \mathcal{E}} \int_{e} \left[\left[\delta p^{n+\gamma}\right]\right] \cdot \left\{\left\{\frac{\grad q}{\rho^{n+\gamma}}\right\}\right\} d\Sigma \nonumber \\
	&+& \sum_{e \in \mathcal{E}}\int_{e} C_{p}\left\{\left\{\frac{1}{\rho^{n+\gamma}}\right\}\right\}_{H}\left[\left[\delta p^{n+\gamma}\right]\right] \cdot \left[\left[q\right]\right]d\Sigma \\
	&=& \sum_{K \in \mathcal{T}_h} \int_{K} \frac{1}{\gamma\Delta t}\mathbf{u}^{n+\gamma,**} \cdot \grad q d\Omega - \sum_{e \in \mathcal{E}} \int_{e} \frac{1}{\gamma\Delta t} \left\{\left\{\mathbf{u}^{n+\gamma,*}\right\}\right\} \cdot \left[\left[q\right]\right]d\Sigma, \nonumber
\end{eqnarray} 
where we set
\begin{equation}
	\sigma^{p}_{e,K} = k^{2}\frac{\text{diam}(e)}{\text{diam}(K)},  
\end{equation}
so that
\begin{equation}
	C_{p} = \frac{1}{2}\left(\sigma^{p}_{e,K^{+}} + \sigma^{p}_{e,K^{-}}\right) \quad \text{if } e \in \mathcal{E}^{I}, \qquad C_{p} = \sigma^{p}_{e,K} \quad \text{if } e \in \mathcal{E}^{B}. 
\end{equation}
The algebraic formulation is once more obtained taking \(q = \psi_{i}, i = 1,\dots,\text{dim}(W_{h})\) and considering the expansion for \(p^{n+\gamma}\) reported in \eqref{eq:fe_representations}. Hence, we get
\begin{equation}
	\left(\frac{M^{2}}{\gamma^{2}\Delta t^{2}}\mathbf{M}_{p}^{n+\gamma} + \mathbf{K}_{p}\right)\mathbf{P}^{n+\gamma} = \mathbf{F}_{p}^{n}.
\end{equation} 
Here, \(\mathbf{P}^{n+\gamma}\) denotes the vector of the degrees of freedom for the pressure. Moreover, we set
\begin{eqnarray}
	\mathbf{M}_{p_{ij}}^{n+\gamma} &=& \sum_{K \in \mathcal{T}_{h}} \int_{K} \psi_{j}\psi_{i} d\Omega \\
	\mathbf{K}_{p_{ij}} &=& \sum_{K \in \mathcal{T}_{h}} \int_{K} \grad\psi_{j} \cdot \grad\psi_{i} d\Omega - \sum_{e \in \mathcal{E}} \int_{e} \left\{\left\{\frac{\grad\psi_{j}}{\rho^{n+\gamma}}\right\}\right\} \cdot \left[\left[\psi_{i}\right]\right]d\Sigma \\
	&-& \sum_{e \in \mathcal{E}} \int_{e} \left[\left[\psi_{j}\right]\right] \cdot \left\{\left\{\frac{\grad\psi_{i}}{\rho^{n+\gamma}}\right\}\right\}d\Sigma + \sum_{e \in \mathcal{E}} \int_{e} C_{p}\left\{\left\{\frac{1}{\rho^{n+\gamma}}\right\}\right\}_{H} \left[\left[\psi_{j}\right]\right] \cdot \left[\left[\psi_{i}\right]\right]d\Sigma \nonumber
\end{eqnarray}
and
\begin{eqnarray}
\mathbf{F}_{p}^{n} &=& \sum_{K \in \mathcal{T}_h} \int_{K} \frac{1}{\gamma\Delta t}\mathbf{u}^{n+\gamma,*} \cdot \grad q d\Omega - \sum_{e \in \mathcal{E}} \int_{e} \frac{1}{\gamma\Delta t} \left\{\left\{\mathbf{u}^{n+\gamma,**}\right\}\right\} \cdot \left[\left[q\right]\right]d\Sigma.
\end{eqnarray}
The second TR-BDF2 stage can be described in a similar manner according to the formulations reported in \eqref{eq:second_stage_levset}, \eqref{eq:second_stage_momentum}, and \eqref{eq:second_stage_helmholtz}.
\\~\\
Finally, we consider the weak formulation for the reinitialization equation for the level set function \eqref{eq:reinitialization_time}:
\begin{eqnarray}
	&&\sum_{K \in \mathcal{T}_{h}}\int_{K} \frac{\phi^{k+1,*}}{\Delta \tau}w d\Omega - \sum_{K \in \mathcal{T}_{h}} \int_{K} u_{c}\phi^{k+1,*}\left(1 - \phi^{k,*}\right)\mathbf{n}_{\Gamma} \cdot \grad w d\Omega \nonumber \\
	&+& \sum_{e \in \mathcal{E}} \int_{e} u_{c}\left\{\left\{\phi^{k+1,*}\left(1 - \phi^{k,*}\right)\mathbf{n}_{\Gamma}\right\}\right\} \cdot \left[\left[w\right]\right]d\Sigma + \sum_{e \in \mathcal{E}} \int_{e} \frac{\tilde{\lambda}^{k}}{2}\left[\left[\phi^{k+1,*}\right]\right] \cdot \left[\left[w\right]\right]d\Sigma \nonumber \\
	&+& \sum_{K \in \mathcal{T}_{h}} \int_{K} u_{c}\beta\varepsilon \left(\grad\phi^{k+1,*} \cdot \mathbf{n}_{\Gamma}\right)\mathbf{n}_{\Gamma} \cdot \grad w d\Omega \\
	&-& \sum_{e \in \mathcal{E}} \int_{e} u_{c}\beta\varepsilon\left\{\left\{\left(\grad\phi^{k+1,*} \cdot \mathbf{n}_{\Gamma}\right)\mathbf{n}_{\Gamma}\right\}\right\} \cdot \left[\left[w\right]\right]d\Sigma \nonumber \\
	&-& \sum_{e \in \mathcal{E}} \int_{e} u_{c}\beta\varepsilon\left\{\left\{\left(\grad v \cdot \mathbf{n}_{\Gamma}\right)\mathbf{n}_{\Gamma}\right\}\right\} \cdot \left[\left[\phi^{k+1,*}\right]\right]d\Sigma + \sum_{e \in \mathcal{E}} \int_{e} C_{\phi}\left[\left[\phi^{k+1,*}\right]\right] \cdot \left[\left[w\right]\right] d\Sigma \nonumber \\
	&=& \sum_{K \in \mathcal{T}_{h}}\int_{K} \frac{\phi^{k,*}}{\Delta \tau}w d\Omega, \nonumber
\end{eqnarray}
where
\begin{equation}
	\tilde{\lambda}^{k} = \max\left(\left|\left(1 - \left(\phi^{k,*}\right)^{+}\right)\mathbf{n}_{\Gamma}^{+} \cdot \mathbf{n}^{+}\right|, \left|\left(1 - \left(\phi^{k,*}\right)^{-}\right)\mathbf{n}_{\Gamma}^{-} \cdot \mathbf{n}^{-}\right|\right).
\end{equation}
Moreover, we set
\begin{equation}
	\sigma^{\phi}_{e,K} = \left(r + 1\right)^{2}\frac{\text{diam}(e)}{\text{diam}(K)},  
\end{equation}
so that
\begin{equation}
	C_{\phi} = \frac{1}{2}\left(\sigma^{\phi}_{e,K^{+}} + \sigma^{\phi}_{e,K^{-}}\right) \quad \text{if } e \in \mathcal{E}^{I}, \qquad C_{\phi} = \sigma^{\phi}_{e,K} \quad \text{if } e \in \mathcal{E}^{B}. 
\end{equation}
One can notice that, following \cite{owkes:2013}, an upwind flux has been employed for the compression term and the SIP has been adopted for the diffusive term. Finally, the algebraic form is obtained considering \(w = \Psi_{i}, i = 1,\dots,\text{dim}(X_{h})\) and the representation in \eqref{eq:fe_representations} so as to obtain
\begin{equation}
	\left(\frac{1}{\Delta\tau}\mathbf{M}_{\phi} + u_{c}\mathbf{C}_{\phi} + \mathbf{A}_{\phi}\right) = \mathbf{F}_{\phi}.
\end{equation}
Here
\begin{eqnarray}
	\mathbf{C}_{\phi_{ij}} &=& -\sum_{K \in \mathcal{T}_{h}} \int_{K} \left(1 - \phi^{k,*}\right)\mathbf{n}_{\Gamma}\psi_{j} \cdot \grad\psi_{i} d\Omega \nonumber \\
	&+& \sum_{e \in \mathcal{E}} \int_{e} \left\{\left\{\left(1 - \phi^{k,*}\right)\mathbf{n}_{\Gamma}\psi_{j}\right\}\right\} : \left[\left[\psi_{i}\right]\right] d\Sigma + \sum_{e \in \mathcal{E}} \int_{e} \frac{\tilde{\lambda}^{k}}{2}\left[\left[\psi_{j}\right]\right] : \left[\left[\psi_{i}\right]\right] d\Sigma \\
	\mathbf{A}_{\phi} &=& \sum_{K \in \mathcal{T}_{h}} \int_{K} u_{c}\beta\varepsilon \left(\grad\psi_{j} \cdot \mathbf{n}_{\Gamma}\right)\mathbf{n}_{\Gamma} \cdot \grad\psi_{i} d\Omega \nonumber \\
	&-& \sum_{e \in \mathcal{E}} \int_{e} u_{c}\beta\varepsilon \left\{\left\{\left(\grad\psi_{j} \cdot \mathbf{n}_{\Gamma}\right)\mathbf{n}_{\Gamma}\right\}\right\} \cdot  \left[\left[\psi_{i}\right]\right] d\Sigma \nonumber \\
	&-& \sum_{e \in \mathcal{E}} \int_{e} u_{c}\beta\varepsilon \left[\left[\psi_{j}\right]\right] \cdot \left\{\left\{\left(\grad\psi_{i} \cdot \mathbf{n}_{\Gamma}\right)\mathbf{n}_{\Gamma}\right\}\right\} d\Sigma \nonumber \\
	&+& \sum_{e \in \mathcal{E}} \int_{e} C_{\phi} \left[\left[\psi_{j}\right]\right] \cdot \left[\left[\psi_{i}\right]\right] d\Sigma.
\end{eqnarray}
Adaptive Mesh Refinement (AMR) will be employed for the numerical tests in Section \ref{sec:tests}. Notice that, AMR does not affect the spatial discretization strategy. More specifically, for faces between cells of different refinement level, the integration is performed from the refined side and a suitable interpolation is performed on the coarse side. Hence, no hanging nodes appear in the implementation of the discrete weak form of the equations.

%%%%%%%%%%%%%%%%%%%%%%%%%%%%%%%%%%%%%%%% Numerical tests %%%%%%%%%%%%%%%%%%%%%%%%
\section{Numerical experiments}
\label{sec:tests} \indent

The numerical method outlined in the previous Sections has been validated in a number of classical test cases for incompressible two-phase flows using the numerical library \texttt{deal.II} \cite{arndt:2022}, whose adaptive mesh refinement capabilities will be employed to enhance resolution close to the interface. The only constraint of the library in the use of non-conforming adaptive meshes is the requirement of not having neighbouring cells with refinement levels differing by more than one. We set \(h = \min\left\{\text{diam}(K) | K \in \mathcal{T}_{h}\right\}\) and we define two Courant numbers, one based on the flow velocity, denoted by \(C_{u}\), and one based on the Mach number, denoted by \(C\):
\begin{equation}
	C_{u} = k\frac{\Delta t U}{h} \qquad C = k\frac{1}{M}\frac{\Delta t}{h},
\end{equation}
where \(U\) is the magnitude of the flow velocity. For the sake of convenience of the reader, we recall here that \(k\) and \(k - 1\) are the polynomial degrees of the finite element spaces chosen for the discretization of velocity and pressure, respectively, whereas \(r\) is the polynomial degree of the finite element space chosen for the discretization of the level set function. We consider \(k = r = 2\) in all the numerical experiments, since the TR-BDF2 is a second order time discretization method, so as to obtain a second order flow solver. Finally, we employ homogenous Neumann boundary conditions for the level set function, while the initial value of the pressure is equal to zero, so as to avoid an extra initial forcing term in the momentum balance.

\subsection{Static bubble}
\label{ssec:static_bubble}

\begin{figure}[h!]
	\centering
	\includegraphics[width = 0.6\textwidth]{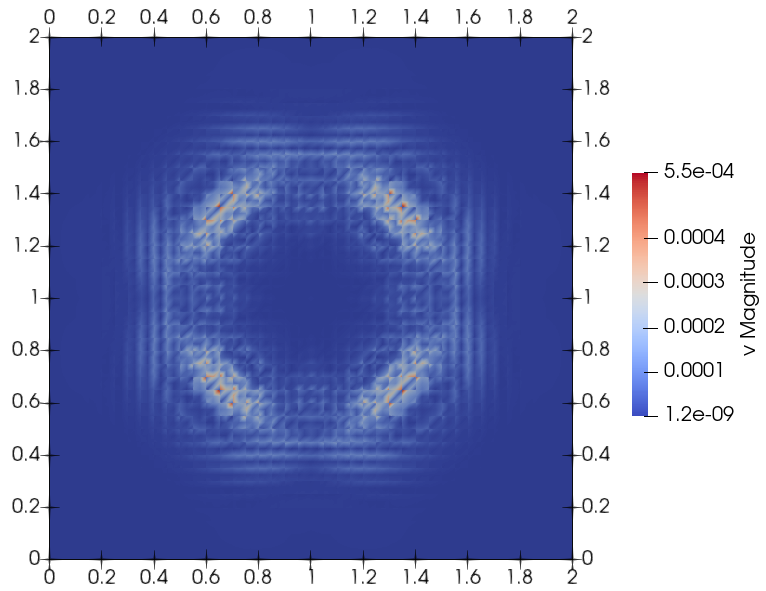}
	\caption{Static bubble test case, velocity magnitude at $t = T_{f} = 250$ using a mesh composed by $40 \times 40$ elements.}
	\label{fig:static_bubble_spurious_currents}
\end{figure}

We consider first the 2D stationary bubble in a zero force field described e.g. in \cite{deshpande:2012, hysing:2006, popinet:1999}. A bubble with \(R = \SI{0.25}{\meter}\) centered in \(\left(x_{0}, y_{0}\right) = \SI[parse-numbers = false]{\left(0.5, 0.5\right)}{\meter}\) in \(\Omega = \SI[parse-numbers=false]{\left(0,1\right)^{2}}{\meter\squared}\) is considered. Following \cite{hysing:2006}, the fluid properties are \(\rho_{1} = \rho_{2} = \SI[parse-numbers=false]{10^{4}}{\kilogram\per\meter\cubed}\), \(\mu_{1} = \mu_{2} = \SI{1}{\kilogram\per\meter\per\second}\), and \(\sigma = \SI{1}{\newton\per\meter}\). Hence, \(\frac{\rho_{1}}{\rho_{2}} = \frac{\mu_{1}}{\mu_{2}} = 1\). In \cite{lafaurie:1994, popinet:1999}, it was conjectured that the amplitude of spurious currents, which arise considering the contribution of the surface tension, is proportional to \(\frac{\sigma}{\mu_{1}}\). Hence, we take \(U_{ref} = \frac{\sigma}{\mu_{1}} = \SI{1}{\meter\per\second}\), so that the Reynolds number equals the Laplace number \(La = \frac{\rho_{1}\sigma L_{ref}}{\mu_{1}^{2}}\). The reference length is equal to \(L_{ref} = 2R = \SI{0.5}{\meter}\), so that the computational domain is \(\Omega = \left(0, 2\right)^{2}\), while the final time is \(T_{f} = 250\), i.e. we consider 250 characteristic times. The artificial speed of sound is set to \(c \approx \SI{1428}{\meter\per\second}\), which is of the same order of magnitude of the speed of sound in the water. Hence, we get \(Re = 5 \times 10^{3}, We = 5 \times 10^{3}\), and \(M = 7 \times 10^{-4}\). We consider no-slip boundary conditions for the velocity and homogeneous Neumann boundary conditions for the pressure. Finally, we set \(\varepsilon = \frac{1}{20}, \Delta\tau = 6.25 \times 10^{-3}, u_{c} = 0.05 u_{max}\) and \(\beta = 0.5\), where \(u_{max}\) is the maximum fluid velocity. The choice to relate \(u_{c}\) with \(u_{max}\) is rather common in the literature, see e.g. \cite{chiu:2011, rusche:2003}. The expected convergence rates are obtained, showing that the numerical method is globally of second order (Table \ref{tab:static_bubble_convergence}). Spurious currents typically appear in the form of vortices around the interface. One can notice that the treatment of the surface tension contribution using the Laplace-Beltrami operator \eqref{eq:laplace_beltrami} reduces the generation of spurious currents, even for a coarse mesh (Figure \ref{fig:static_bubble_spurious_currents}).

\begin{table}[!h]
	\centering
	\footnotesize
	\begin{tabularx}{\textwidth}{rrrcrcrc}
		\toprule
		$N_{el}$ & $\Delta t$ & $\left\|\mathbf{u}\right\|_{H^{1}(\Omega)}$ & Conv. rate $H^{1}$ & $\left\|\mathbf{u}\right\|_{L^{2}(\Omega)}$ & Conv. rate $L^{2}$ & $\left\|\mathbf{u}\right\|_{L^{\infty}(\Omega)}$ & Conv. rate $L^{\infty}$ \\
		\midrule
		$20$ & $1$ & $4.37 \times 10^{-2}$ & & $4.81 \times 10^{-3}$ & & $4.77 \times 10^{-3}$ & \\ 
		\midrule
		$40$ & $0.5$ & $7.59 \times 10^{-3}$ & $2.5$ & $7.48 \times 10^{-5}$ & $2.7$ & $2.81 \times 10^{-4}$ & $ 4.1$ \\
		\midrule
		$80$ & $0.25$ & $1.14 \times 10^{-3}$ & $2.7$ & $7.23 \times 10^{-6}$ & $3.4$ & $3.01 \times 10^{-5}$ & $ 3.2$ \\
		\midrule
		$160$ & $0.125$ & $2.43 \times 10^{-4}$ & $\mathbf{2.2}$ & $8.64 \times 10^{-7}$ & $\mathbf{3.1}$ & $3.32 \times 10^{-6}$ & $\mathbf{3.2}$ \\
		\bottomrule
	\end{tabularx}
	\caption{Static bubble test case, convergence rates at fixed $\frac{U_{ref}\Delta t}{h}$ in norm $H^{1}$, $L^{2}$, and $L^{\infty}$. $N_{el}$ denotes the number of elements along each direction.}
	\label{tab:static_bubble_convergence}
\end{table}

\subsection{Rayleigh-Taylor instability}
\label{ssec:RT}

The Rayleigh-Taylor instability is a well known test case in which an heavier fluid penetrates a lighter fluid under the action of gravity. We consider the configuration presented e.g. in \cite{bassi:2020, haghshenas:2017, popinet:1999}, for which \(\rho_{1} = \SI{1.225}{\kilogram\per\meter\cubed}\) and \(\rho_{2} = \SI{0.1694}{\kilogram\per\meter\cubed}\), corresponding to the density of air and helium, respectively, whereas \(\mu_{1} = \mu_{2} = \SI{0.00313}{\kilogram\per\meter\per\second}\). The effect of surface tension is neglected. Following \cite{tryggvason:1988}, we consider as reference length, the computational width of the box \(W\) and as reference time, the time scale of wave growth, equal to \(t_{ref} = \sqrt{\frac{W}{Ag}}\), where \(g = \SI{9.81}{\meter\per\second\squared}\) and \(A = \frac{\rho_{1} - \rho_{2}}{\rho_{1} + \rho_{2}}\) is the Atwood number. Hence, we obtain the following relations:
\begin{equation}
	U_{ref} = \sqrt{AgW} \qquad Re = \frac{\rho_{1}\sqrt{AgW}W}{\mu_{1}} \qquad Fr = \sqrt{A}. 
\end{equation}
We consider \(W = \SI{1}{\meter}\), so as to obtain a computational domain \(\Omega = \left(0,1\right) \times \left(0,4\right)\). Hence, we get \(A \approx 0.757, t_{ref} \approx \SI{0.367}{\second}, U_{ref} \approx \SI{2.725}{\meter\per\second}, Re \approx 1066.55,\) and \(Fr \approx 0.87\). We take \(M = 0.008\), corresponding to \(c \approx \SI{343}{\meter\per\second}\), namely the speed of sound in air. The final time is \(T_{f} = 2.45\). No-slip boundary conditions are prescribed on top and bottom walls, whereas periodic boundary conditions are imposed along the horizontal direction. The pressure is prescribed to be zero on the upper wall, while homogeneous Neumann boundary conditions are employed for the bottom wall. The initial velocity field is zero, whereas the initial level set function is
\begin{equation}
	\phi(0) = \frac{1}{1 + \exp\left(\frac{2 + 0.05\cos(2\pi x) - y}{\varepsilon}\right)}.
\end{equation}
The computational grid is composed by \(160 \times 640\) elements, whereas the time step is \(\Delta t \approx 1.63 \times 10^{-3}\), yielding a maximum advective Courant number \(C \approx 1.36\) and an acoustic Courant number \(C \approx 65.3.\) Finally, we set \(\varepsilon = h = \frac{1}{160}, \Delta\tau = 0.05h, u_{c} = 0.0125u_{max}\), and \(\beta = 1\). Figure \ref{fig:RT_At0,76_phi} shows the development of the interface at \(t = T_{f}\), where one can easily notice the expected main behaviour of the Rayleigh-Taylor instability: as the heavier fluid penetrates the lighter one, the interface begins to roll up along the sides of the spike giving the typical ``mushroom'' shape. Obtained results are similar to those in literature, see e.g. \cite{haghshenas:2017, popinet:1999, puckett:1997, sheu:2009}. Moreover, for the sake of completeness, we report in Figure \ref{fig:RT_At0,76_area} the evolution of the relative variation of the area for the lighter fluid, defined as
\begin{equation}
	\frac{\left|\Omega_{2}(t) - \Omega_{2}(0)\right|}{\Omega_{2}(0)}.
\end{equation}
The maximum relative variation is 0.034 \%, showing that CLS method preserves the area quite well. The results compare well with those in \cite{salih:2006}, in which a loss of around $1.1\%$ was experienced.

\begin{figure}[!h]
	\centering
	\includegraphics[width=0.45\textwidth]{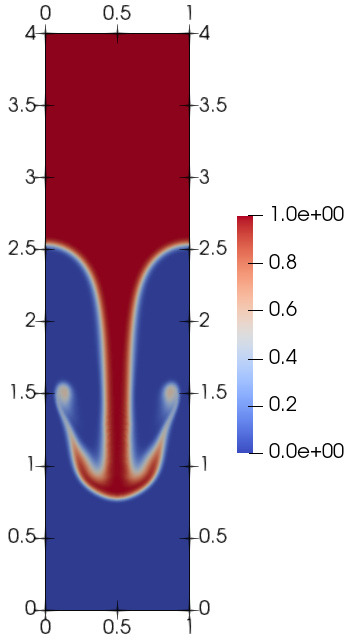}
	\caption{Rayleigh-Taylor instability, contour plot of the level set function at $t = T_{f} = 2.45$.}
	\label{fig:RT_At0,76_phi}
\end{figure}

\begin{figure}[!h]
	\centering
	\includegraphics[width=0.75\textwidth]{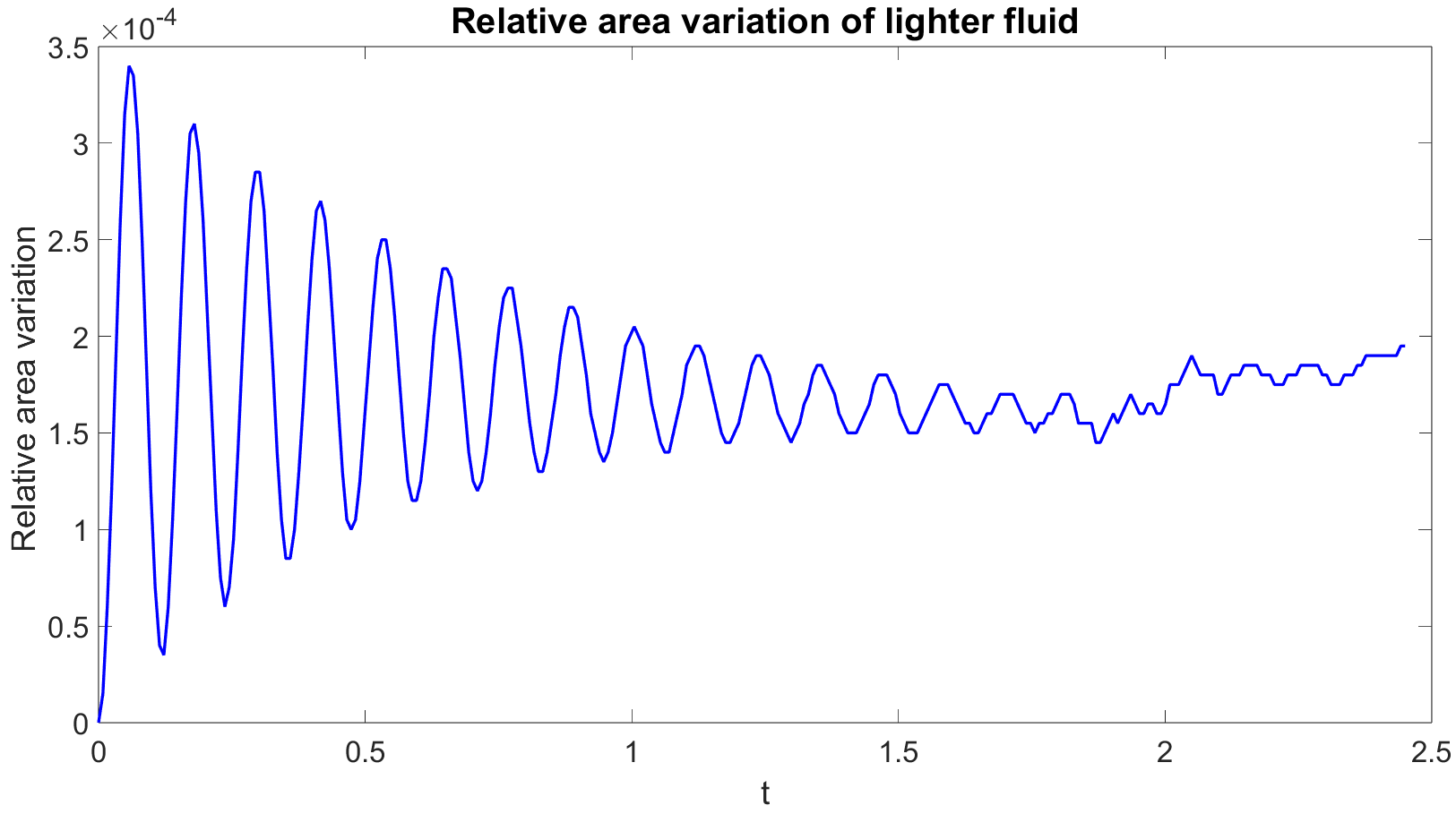}
	\caption{Rayleigh-Taylor instability, evolution of the relative variation of the area for the lighter fluid.}
	\label{fig:RT_At0,76_area}
\end{figure}

\begin{figure}[!h]
	\centering
	\includegraphics[width=0.25\textwidth]{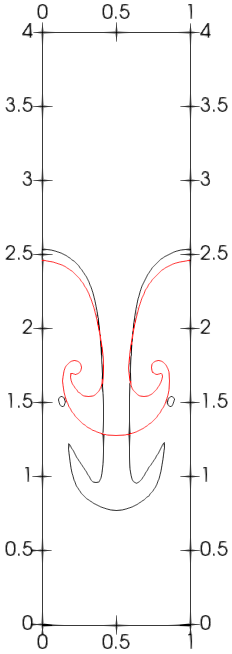}
	\caption{Rayleigh-Taylor instability, comparison between $A \approx 0.757$ and $A \approx 0.5$. The black line shows the interface for $A \approx 0.757$, whereas the red line refers to the interface for $A \approx 0.5$.}
	\label{fig:RT_phi_comparison}
\end{figure}

An interesting analysis regards the influence of the Atwood number. We fix \(\rho_{2} = \SI{0.408}{\kilogram\per\meter\cubed}\), so as to obtain \(A \approx 0.5\). As a consequence, we obtain \(t_{ref} \approx \SI{0.451}{\second}, U_{ref} \approx \SI{2.215}{\meter\per\second}, Re \approx 867.05, Fr \approx 0.71\), and \(M = 0.006\). We set the final time \(T_{f} = 2\), so that the same final dimensional time of the previous configuration is achieved. The chosen time step is \(\Delta t = 2.5 \cdot 10^{-3}\). One can easily notice from Figure \ref{fig:RT_phi_comparison} that, with higher Atwood number, the roll up effect is enhanced. This points out the earlier appearance of the Kelvin-Helmholtz instability, due to the development of short wavelength perturbations along the fluid interface.

\begin{figure}[!h]
	\begin{subfigure}{0.5\textwidth}
		\centering
		\includegraphics[width=0.94\textwidth]{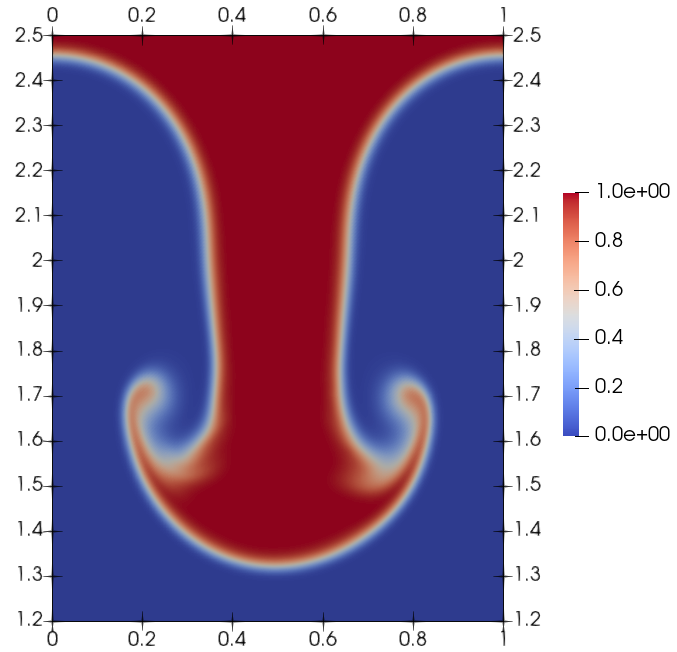}
	\end{subfigure}
	\begin{subfigure}{0.5\textwidth}
		\centering
		\includegraphics[width=0.75\textwidth]{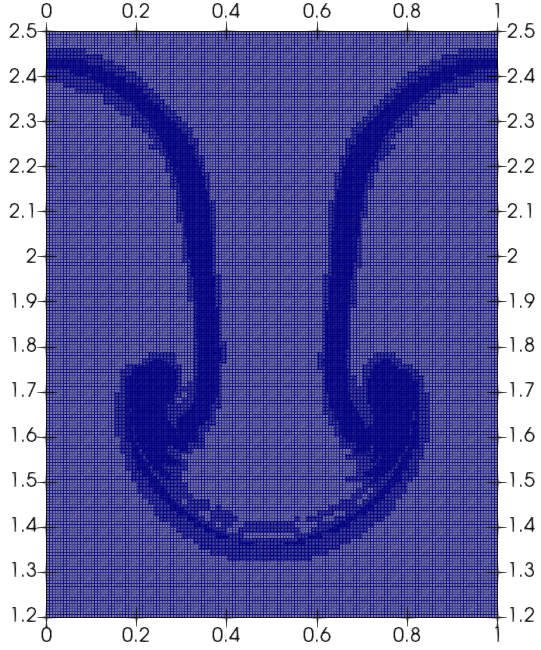}
	\end{subfigure}
	\caption{Rayleigh-Taylor instability at $A \approx 0.5$. Left: contour plot of the level set function at $t = T_{f} = 2$. Right: computational grid at $t = T_{f} = 2$.}
	\label{fig:RT_adaptive}
\end{figure}

\begin{figure}[!h]
	\begin{subfigure}{0.5\textwidth}
		\centering
		\includegraphics[width=0.8\textwidth]{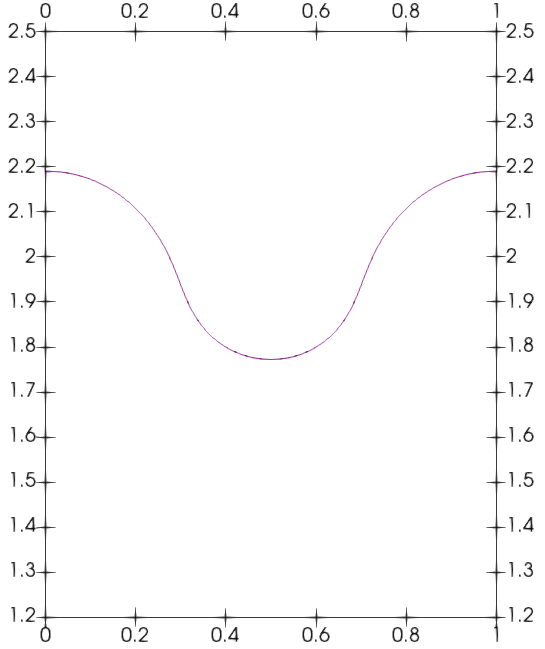}
	\end{subfigure}
	\begin{subfigure}{0.5\textwidth}
		\centering
		\includegraphics[width=0.8\textwidth]{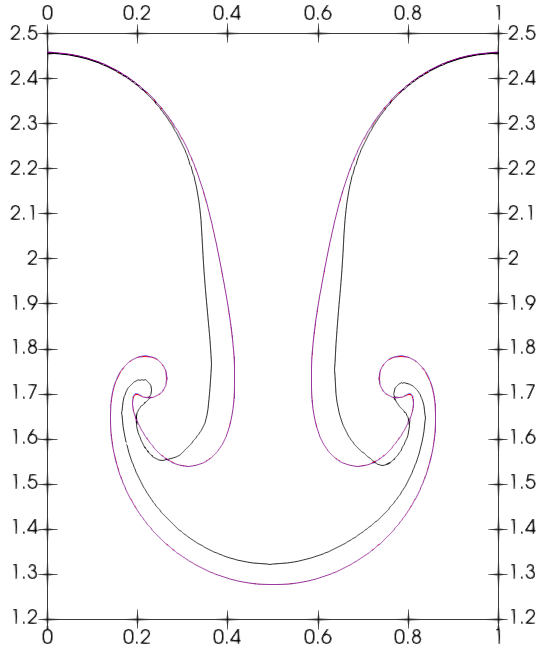}
	\end{subfigure}
	\caption{Rayleigh-Taylor instability, comparison at $A \approx 0.5$ between the fixed grid simulation and the adaptive grid one. Left: results at $t = \frac{T_{f}}{2} = 1$. Right: results at $t = T_{f} = 2$. The black lines show the interface obtained with the adaptive grid, the red lines refer to the interface obtained with the fixed grid composed by $160 \times 640$ elements, whereas the blue lines report the results obtained with a fixed grid composed by $320 \times 1280$ elements.}
	\label{fig:RT_phi_comparison_adaptive}
\end{figure}

We employ now the \(h\)-adaptive version of the scheme for the latter configuration. More specifically, we define for each element \(K\) the quantity
\begin{equation}\label{eq:eta_K}
	\eta_{K} = \max_{i \in \mathcal{N}_{K}}\left|\grad\phi\right|_{i},
\end{equation}
which acts as local refinement indicator. Here \(\mathcal{N}_{K}\) denotes the set of nodes over the element \(K\). We allow to refine when \(\eta_{K}\) exceeds \(10\) and to coarsen below \(5\). The initial grid is composed by \(80 \times 320\) elements and we allow up to two local refinements, so as to obtain \(h = \frac{1}{320}\) and a maximum resolution which would correspond to a \(320 \times 1280\) uniform grid. As one can notice from Figure \ref{fig:RT_adaptive}, the refinement criterion is able to increase the resolution only in correspondence of the interface between the two fluids. The final grid consists of \(42850\) elements, corresponding to around 40 \% of elements of the fixed uniform grid. A solution using the full resolution uniform grid, i.e. a \(320 \times 1280\) grid, has been also computed for further comparison. One can easily notice that at \(t = \frac{T_{f}}{2}\) the interfaces are indistinguishable, whereas at \(t = T_{f}\) a slightly different development of the instability appears between the simulation with the adaptive grid and those with the two fixed grids (Figure \ref{fig:RT_phi_comparison_adaptive}). Moreover, one can notice that the same profiles of the interface are obtained for the two fixed grids, meaning that we have achieved grid independence. Since we are analyzing a fluid mechanic instability, every small variation in the flow corresponds to large variations, and, therefore, it is difficult to say which solution is the more reliable. Similar results and considerations have been reported for a Kelvin-Helmholtz instability in \cite{orlando:2023a}. Finally, for what concerns the computation efficiency, AMR leads to a computational time saving of around 90 \% for a given maximum spatial resolution (Table \ref{tab:RT_times}).

\begin{table}[h!]
	\centering
	\footnotesize
	\begin{tabularx}{0.75\columnwidth}{rc}
		\toprule
		Grid & WT$[\SI{}{\second}]$ \\
		\midrule
		$160 \times 640$ (fixed) & 2110 \\
		\midrule
		$320 \times 1280$ (fixed) & 14400 \\
		\midrule
		adaptive ($N_{el}$ = 42850 at $t = T_{f}$, maximum resolution $320 \times 1280$) & 1830 \\
		\bottomrule
	\end{tabularx}
	\caption{Rayleigh-Taylor instability, wall-clock times for the different configurations at $A \approx 0.5$.}
	\label{tab:RT_times}
\end{table}

\subsection{Rising bubble benchmark}
\label{ssec:rising_bubble}

The rising bubble benchmark is a well-established test case for the validation of numerical methods for incompressible two-phase flows \cite{hysing:2009}. More specifically, the evolution of the shape, position and velocity of the center of mass of a rising bubble is compared against the reference solution in \cite{hysing:2009}. Two configurations are considered with the corresponding physical parameters and non-dimensional numbers listed in Table \ref{tab:rising_bubble_param} and \ref{tab:rising_bubble_adim_param}, respectively. The bubble occupies the subdomain \(\Omega_{2}\). Following \cite{hysing:2009}, we set \(L_{ref} = 2r_{0} = \SI{0.5}{\meter}\) and \(U_{ref} = \sqrt{g L_{ref}} = \SI{0.7}{\meter\per\second}\). We consider as domain \(\Omega = \left(0, L_{x}\right) \times \left(0, L_{y}\right)\), with \(L_{x} = 2\) and \(L_{y} = 4\), whereas the final time is \(T_{f} = 4.2.\) No-slip boundary conditions are imposed on the top and bottom boundaries, with homogeneous Neumann boundary conditions for the pressure. Periodic conditions are prescribed in the horizontal direction. The initial velocity field is zero. Finally, the initial level set function is described by the following relation:
\begin{equation}
	\phi(0) = \frac{1}{1 + \exp\left(\frac{R - \sqrt{\left(x - x_{0}\right)^{2} + \left(y - y_{0}\right)^{2}}}{\varepsilon}\right)},
\end{equation}
with \(R = 1, x_{0} = y_{0} = 1\). We compute as reference quantities the position  \(\mathbf{x}_{c}\), the velocity \(\mathbf{u}_{c}\) of the center of mass, and the so-called degree of circularity \(\chi\), defined respectively as
\begin{eqnarray}
	\mathbf{x}_{c} &=& \frac{\int_{\Omega_{2}}\mathbf{x}d\Omega}{\int_{\Omega_{2}}d\Omega} = \frac{\int_{\Omega_{2}}\mathbf{x}d\Omega}{\left|\Omega_{2}\right|} \\
	\mathbf{u}_{c} &=& \frac{\int_{\Omega_{2}}\mathbf{u}d\Omega}{\int_{\Omega_{2}}d\Omega} = \frac{\int_{\Omega_{2}}\mathbf{u}d\Omega}{\left|\Omega_{2}\right|} \\
	\chi &=& \frac{2\sqrt{\pi\left|\Omega_{2}\right|}}{P_b},
\end{eqnarray}
where \(\Omega_{2}\) is the subdomain occupied by the bubble, \(\left|\Omega_{2}\right|\) is the area of the bubble, and \(P_{b}\) is its perimeter. The degree of circularity is the ratio between the perimeter of a circle with the same area of the bubble and the current perimeter of the bubble itself. For a perfectly circular bubble, the degree of circularity is equal to one and then decreases as the bubble deforms itself. Since \(\phi\) is a regularized Heaviside function, we can compute the reference quantities as follows:
\begin{eqnarray}
	\mathbf{x}_{c} &\approx& \frac{\int_{\Omega}\mathbf{x}\left(1- \phi\right)d\Omega}{\int_{\Omega}\left(1- \phi\right)d\Omega} \\
	\mathbf{u}_{c} &\approx& \frac{\int_{\Omega}\mathbf{u}\left(1- \phi\right)d\Omega}{\int_{\Omega}\left(1- \phi\right)d\Omega} \\
	\chi &\approx& \frac{2\sqrt{\pi \int_{\Omega}\left(1- \phi\right)d\Omega}}{\int_{\Omega}\left|\grad\phi\right|d\Omega}.
\end{eqnarray} 
\begin{table}[!h]
	\centering
	\begin{tabularx}{\textwidth}{rcccccc}
		\toprule
		Test case & \(\rho_{1} \left[\SI{}{\kilogram\per\meter\cubed}\right]\) & \(\rho_{2} \left[\SI{}{\kilogram\per\meter\cubed}\right]\) & \(\mu_{1} \left[\SI{}{\kilogram\per\meter\per\second}\right]\) & \(\mu_{2} \left[\SI{}{\kilogram\per\meter\per\second}\right]\) & \(g \left[\SI{}{\meter\per\second\squared}\right]\) & \(\sigma \left[\SI{}{\kilogram\per\second\squared}\right]\) \\
		\midrule
		Config. 1 & \(1000\) & \(100\) & \(10\) & \(1\) & \(0.98\) & \(24.5\) \\ 
		\midrule
		Config. 2 & \(1000\) & \(1\) & \(10\) & \(0.1\) & \(0.98\) & \(1.96\) \\
		\bottomrule
	\end{tabularx}
	\caption{Physical parameters defining the configurations from rising bubble test case (data from \cite{hysing:2009}).}
	\label{tab:rising_bubble_param}
\end{table}

\begin{table}[!h]
	\centering
	\begin{tabularx}{0.5\textwidth}{rccccc} 
		\toprule
		Test case & \(Re\) & \(Fr\) & \(We\) & \(\rho_{2}/\rho_{1}\) & \(\mu_{2}/\mu_{1}\) \\
		\midrule
		Config. 1 & \(35\) & \(1\) & \(10\) & \(10^{-1}\) & \(10^{-1}\) \\ 
		\midrule
		Config. 2 & \(35\) & \(1\) & \(125\) & \(10^{-3}\) & \(10^{-2}\) \\
		\bottomrule
	\end{tabularx}	
	\caption{Non-dimensional numbers defining the configurations from rising bubble test case (data from \cite{hysing:2009}).}
	\label{tab:rising_bubble_adim_param}
\end{table}
We start with the first configuration and we set \(M = 0.0005\), corresponding to \(c = \SI{1400}{\meter\per\second}\), which is of the order of magnitude of the speed of sound in water. The computational grid is composed by \(320 \times 640\) elements, leading to \(h = \frac{1}{160}\), whereas the time step is \(\Delta t = 6 \cdot 10^{-3}\), yielding a maximum advective Courant number \(C_{u} \approx 1.4\) and an acoustic Courant number \(C = 1920\). Finally, we set \(\varepsilon = h, \Delta\tau = 0.05h, u_{c} = 0.05 u_{max}\) and \(\beta = 0.5\). We point out here the fact that results in the Figures have been compared with the results of Group 2 in \cite{hysing:2009}. One can easily notice that we are able to recover the reference shape of the bubble at \(t = T_{f}\) (Figure \ref{fig:rising_bubble_case1_shape}). A good qualitative agreement is established for the evolution of the degree of circularity, with only slightly lower values for our numerical results (Figure \ref{fig:rising_bubble_case1}). {A maximum relative error of the order of $10^{-1}$ is established with respect to the results of Group 2 in \cite{hysing:2009}. The center of mass reaches \(y_{c} = 2.156\), which is in good agreement with the value \(y_{c} = 2.162 \pm 0.002\) reported in \cite{hysing:2009} and implies a relative error of the order of $10^{-3}$ (Figure \ref{fig:rising_bubble_case1}). Finally, the maximum rise velocity of the center of mass is \(v_{c} = 0.3461\), which is again in good agreement with the value \(v_{c} = 0.3456 \pm 0.0003\) present in \cite{hysing:2009} and implies again a relative error of the order of $10^{-3}$ (Figure \ref{fig:rising_bubble_case1}). 

\begin{figure}[!h]
	\begin{subfigure}{0.5\textwidth}
		\centering
		\includegraphics[width=0.9\textwidth]{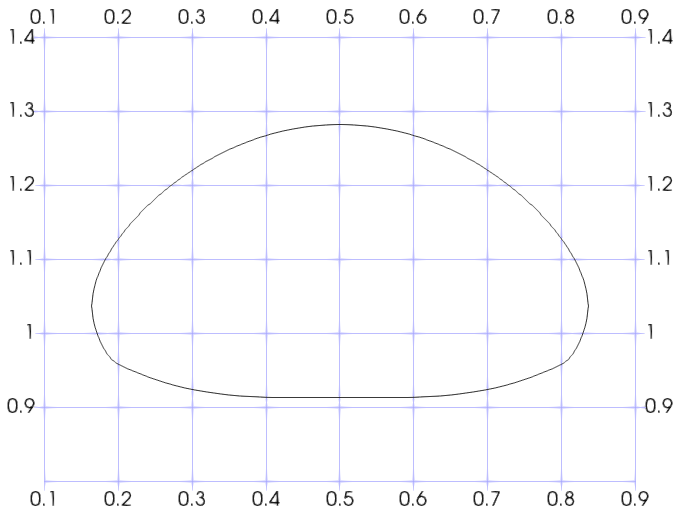}
	\end{subfigure}
	\begin{subfigure}{0.5\textwidth}
		\centering
		\includegraphics[width=0.9\textwidth]{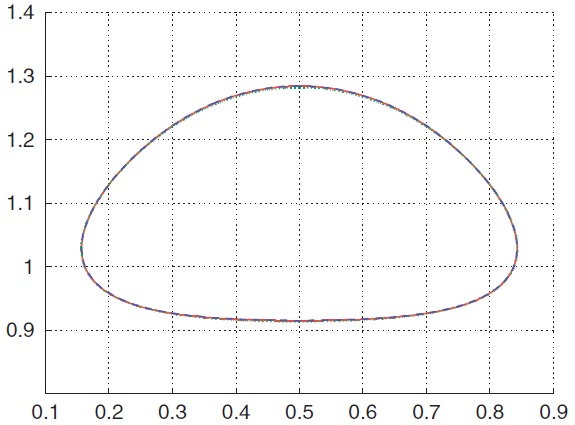}
	\end{subfigure}
	\caption{Rising bubble benchmark, configuration 1, shape of bubble at $t = T_{f} = 4.2$. Left: numerical simulation. Right: image from \cite{hysing:2009}. Bounds have been rescaled by $L_{ref}$ for the sake of comparison with reference results.}
	\label{fig:rising_bubble_case1_shape}
\end{figure}

\begin{figure}[!h]
	\begin{subfigure}{0.5\textwidth}
		\centering
		\includegraphics[width=0.9\textwidth]{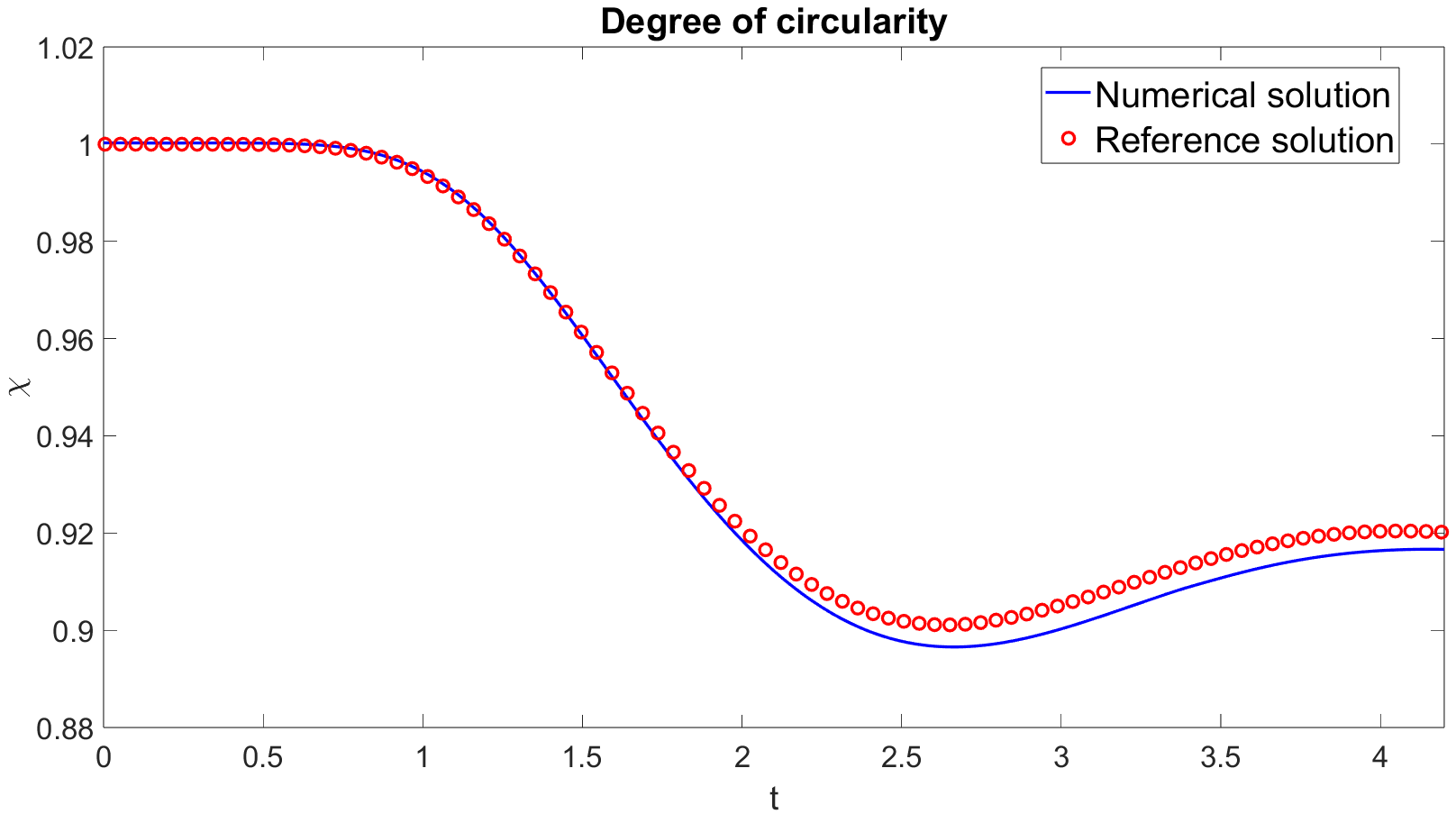} a)
	\end{subfigure}
	\begin{subfigure}{0.5\textwidth}
		\centering
		\includegraphics[width=0.9\textwidth]{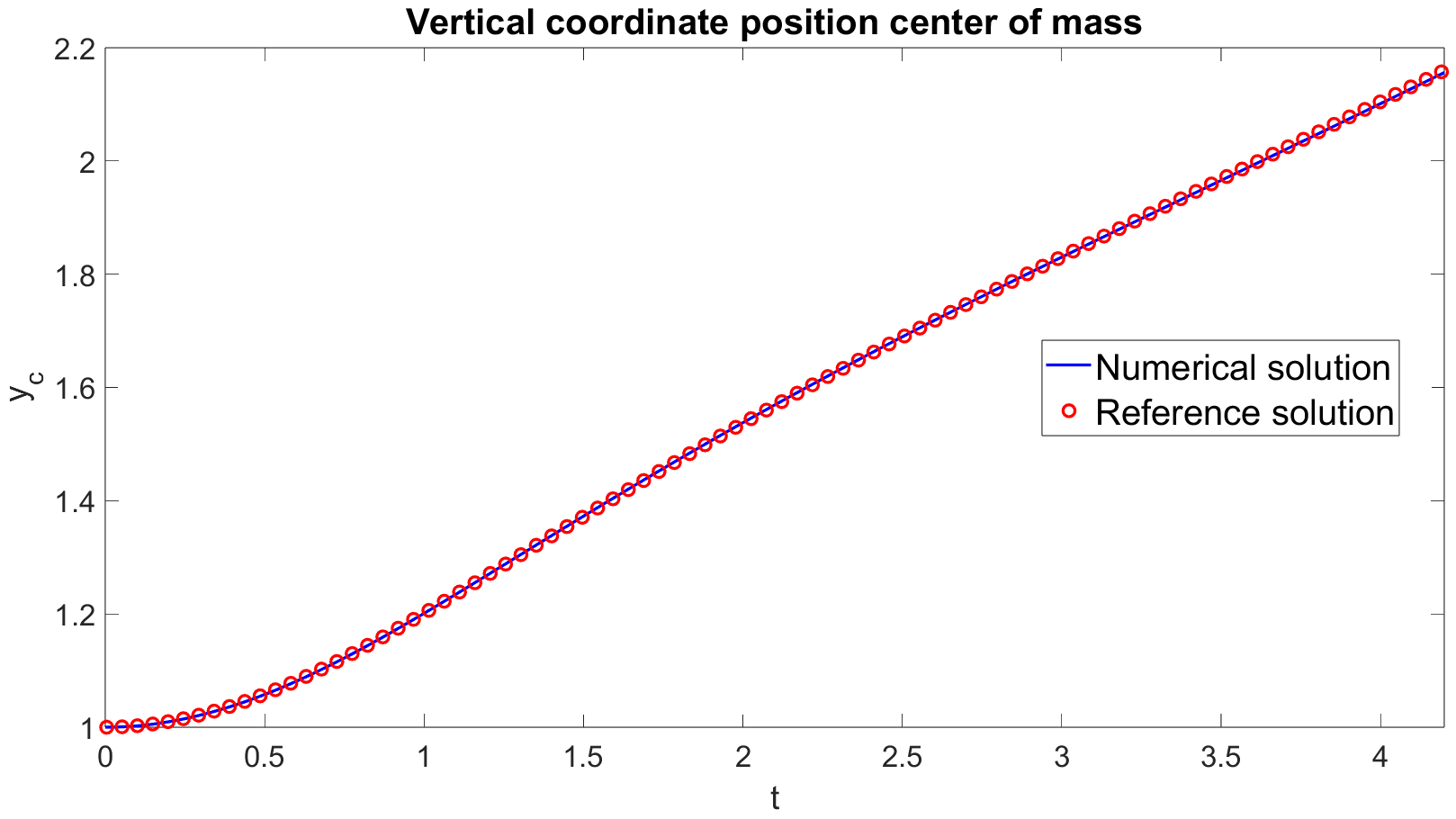} b)
	\end{subfigure}
	\begin{subfigure}{0.5\textwidth}
		\centering
		\includegraphics[width=0.9\textwidth]{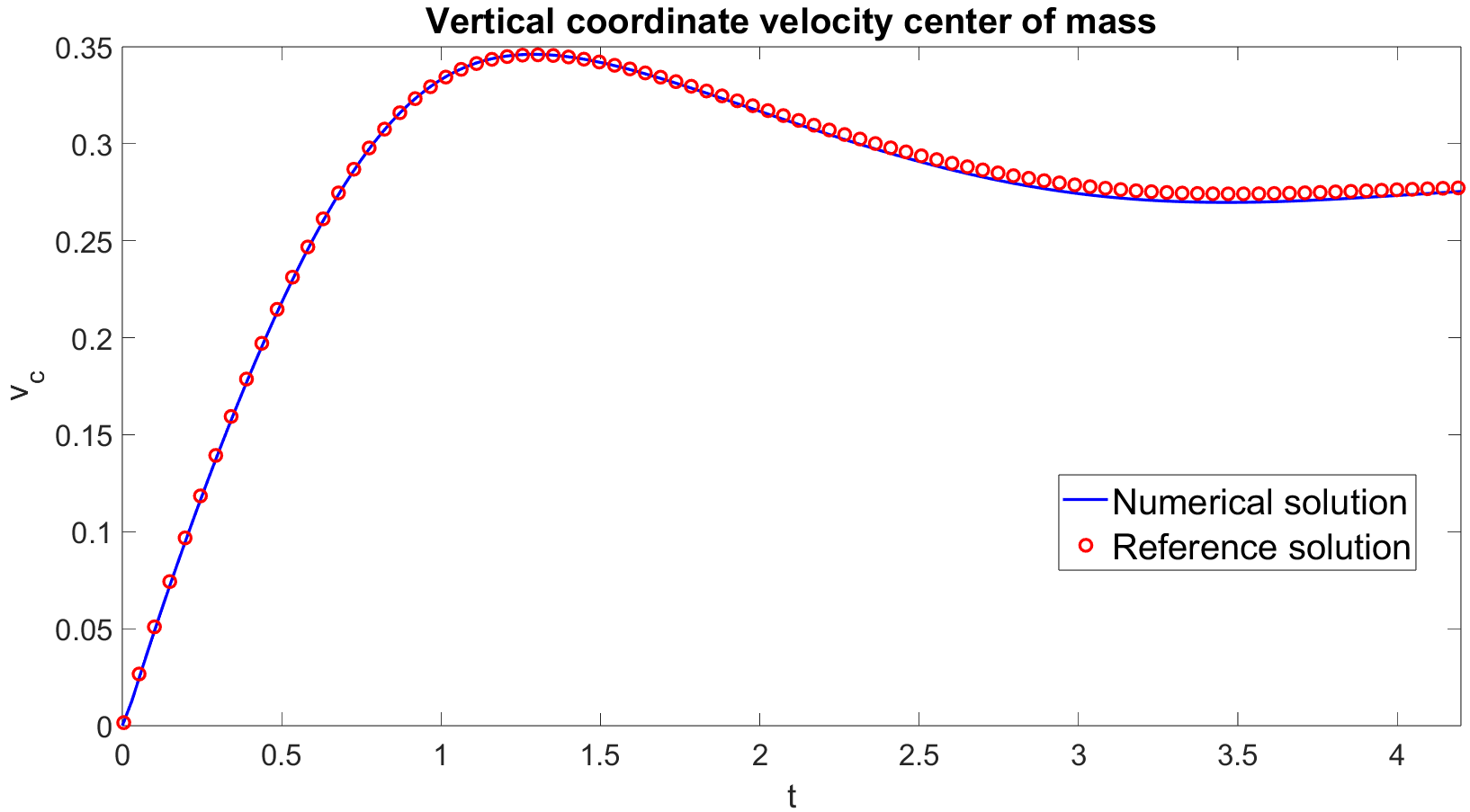} c)
	\end{subfigure}	
	\caption{Rising bubble benchmark, configuration 1, a) degree of circularity, b) vertical coordinate position of the center of mass, c) vertical coordinate velocity of the center of mass. The blue lines denote the results obtained with our method, while the red dots show the reference results from \cite{hysing:2009}. Reference data for position and velocity have been rescaled by $L_{ref}$ and $U_{ref}$, respectively, for the sake of comparison.}
	\label{fig:rising_bubble_case1}
\end{figure}

We analyze now the second configuration. The time step is \(\Delta t = 5 \cdot 10^{-3}\), yielding a maximum advective Courant number \(C_{u} \approx 1.4\) and an acoustic Courant number \(C = 1600\). We also set \(\varepsilon = h = \frac{1}{160}, \Delta\tau = 0.05h = 3.125 \times 10^{-4}, u_{c} = 0.0125 u_{max}\) and \(\beta = 2\). The bubble develops a non-convex shape with thin filaments (Figure \ref{fig:rising_bubble_case2_shape}). The solutions given in \cite{hysing:2009} are different and, in some cases, the thin filaments tend to break off, although it is unclear if such a phenomenon should be observed in the current two-dimensional setting. The obtained profile is however in good agreement with that of Group 2 in \cite{hysing:2009}. Figure \ref{fig:rising_bubble_case2} shows the evolution of the degree of circularity, of the vertical coordinate of the position of the center of mass, and of the vertical coordinate of the velocity of the center of mass. A good qualitative agreement is established for the quantities of interest, even though deviations from the chosen reference solution are visible. In particular, differences appear for the degree of circularity starting from \(t \approx 2.5\), when the thin filaments start developing. Moreover, the second peak for the rising velocity reaches a lower value. As mentioned above, there is no clear agreement concerning the thin filamentary regions, and, therefore, their development can strongly affect computations of the reference quantities and can lead to different numerical results.

\begin{figure}[!h]
	\begin{subfigure}{0.5\textwidth}
		\centering
		\includegraphics[width=0.9\textwidth]{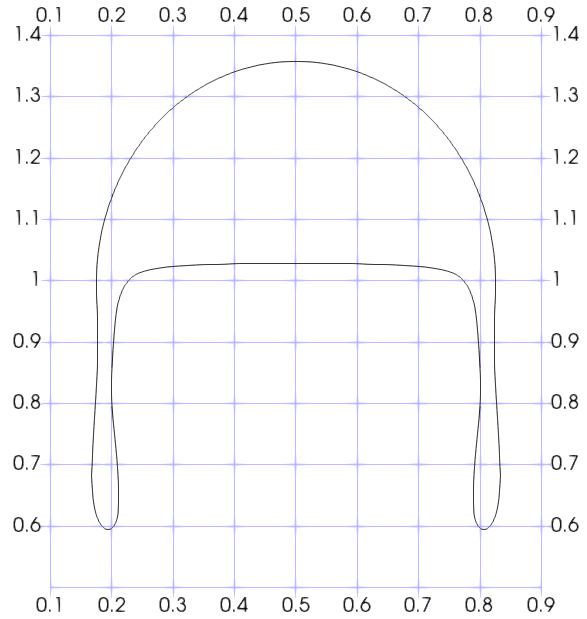}
	\end{subfigure}
	\begin{subfigure}{0.5\textwidth}
		\centering
		\includegraphics[width=0.88\textwidth]{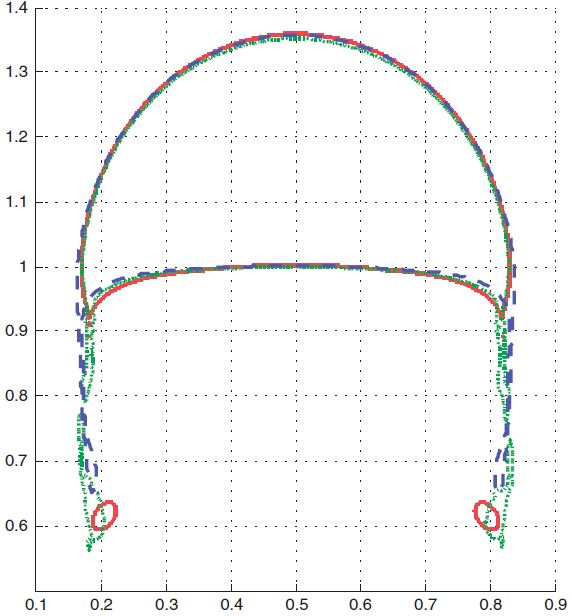}
	\end{subfigure}
	\caption{Rising bubble benchmark, configuration 2, shape of bubble at $t = T_{f} = 4.2$. Left: numerical simulation. Right: image from \cite{hysing:2009}. Bounds have been rescaled by $L_{ref}$ for the sake of comparison with the reference results.}
	\label{fig:rising_bubble_case2_shape}
\end{figure}

\begin{figure}[!h]
	\begin{subfigure}{0.5\textwidth}
		\centering
		\includegraphics[width=0.9\textwidth]{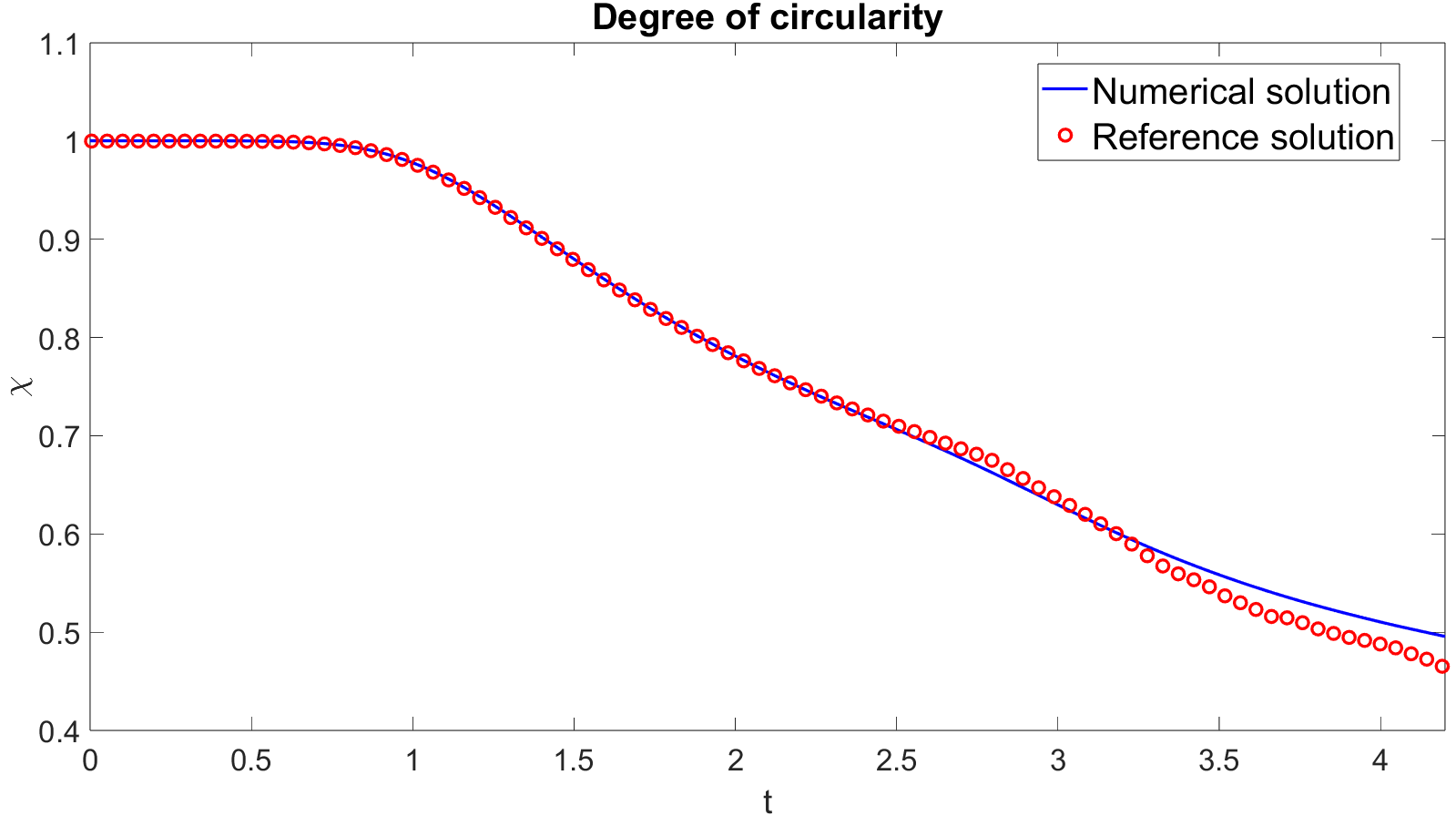} a)
	\end{subfigure}
	\begin{subfigure}{0.5\textwidth}
		\centering
		\includegraphics[width=0.9\textwidth]{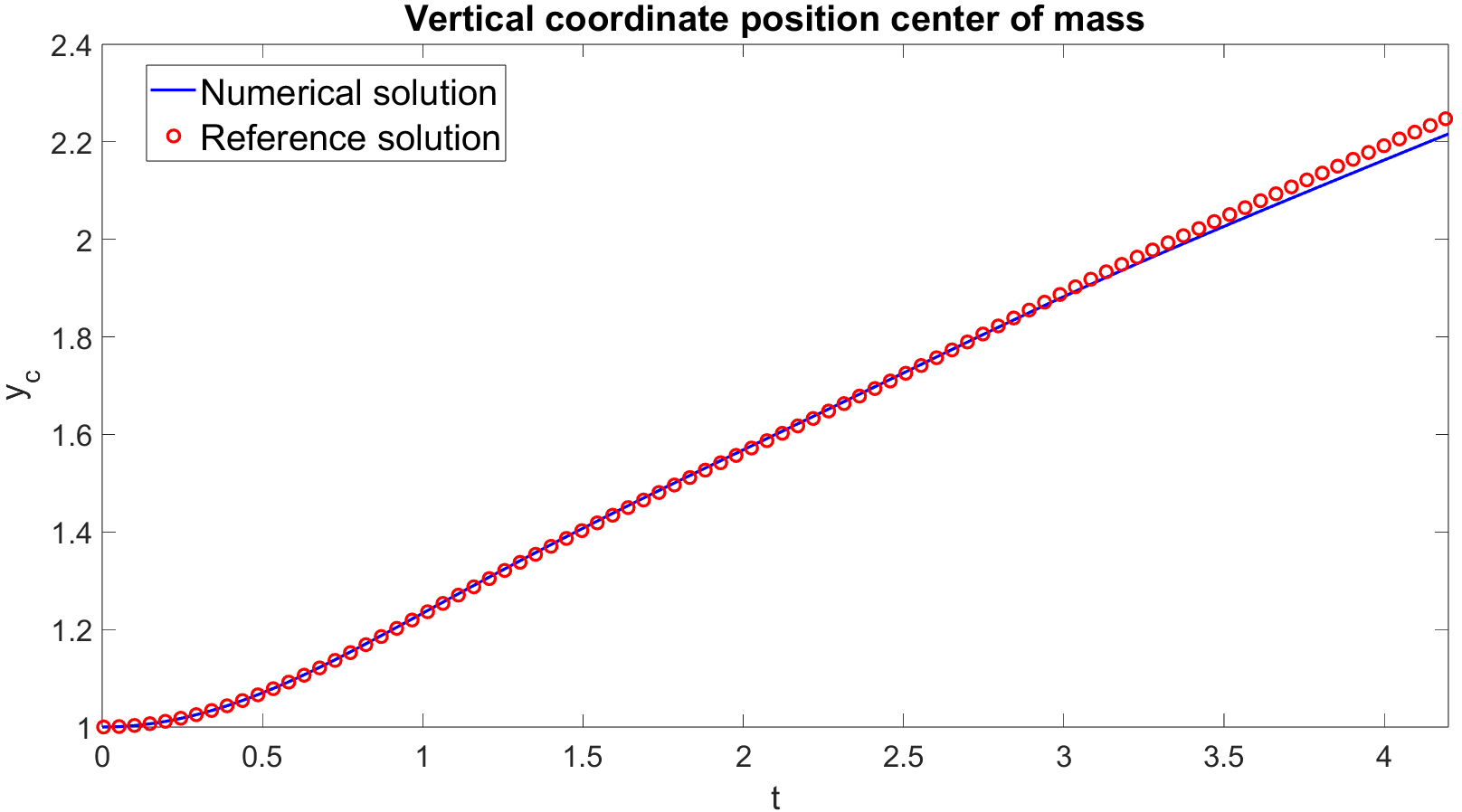} b)
	\end{subfigure}
	\begin{subfigure}{0.5\textwidth}
		\centering
		\includegraphics[width=0.9\textwidth]{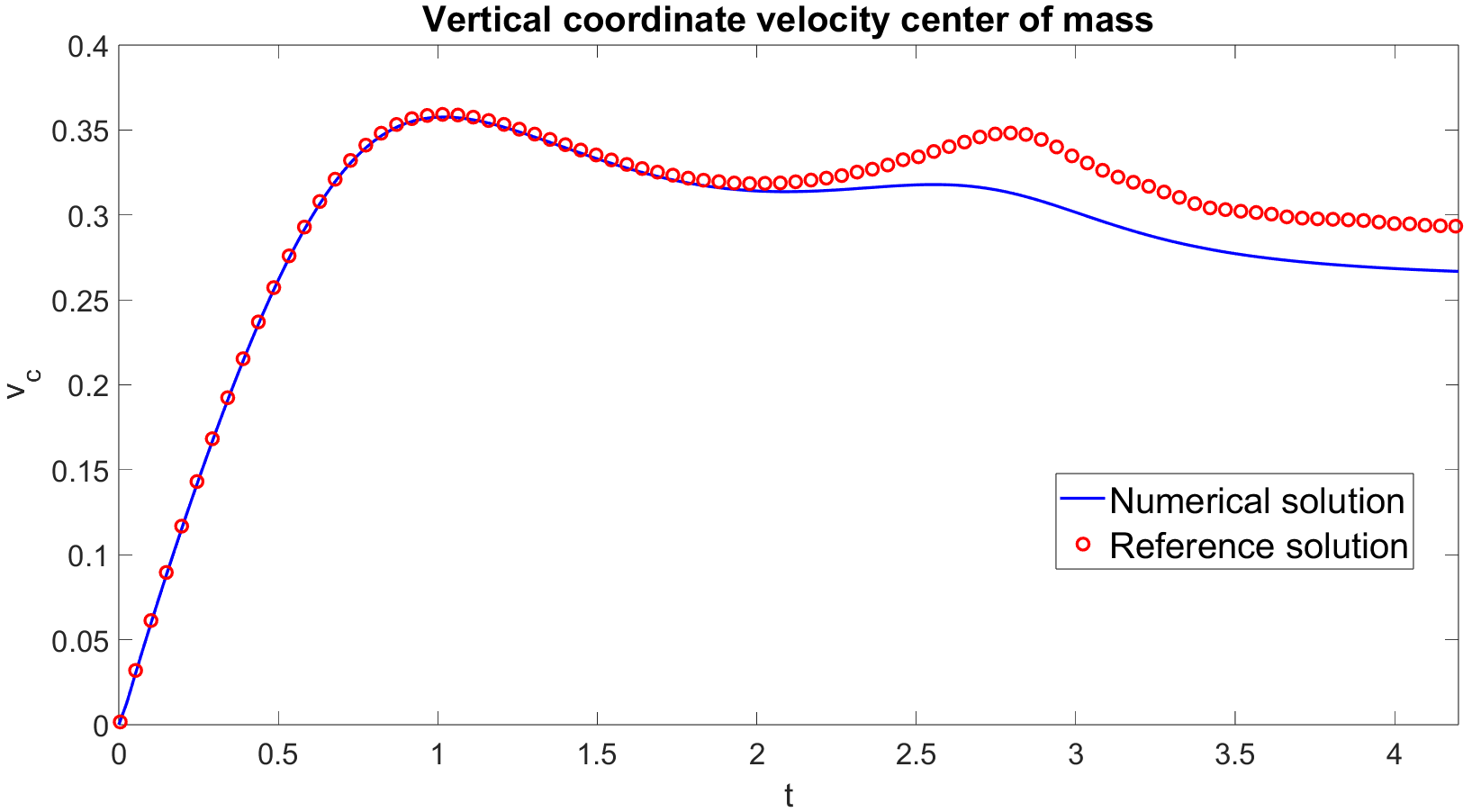} c)
	\end{subfigure}
	\caption{Rising bubble benchmark, configuration 2, a) degree of circularity, b) vertical coordinate position of the center of mass, c) vertical coordinate velocity of the center of mass. The blue lines denote the results obtained with our method, while the red dots show the reference results from \cite{hysing:2009}. Reference data for position and velocity have been rescaled by $L_{ref}$ and $U_{ref}$, respectively, for the sake of comparison.}
	\label{fig:rising_bubble_case2}
\end{figure}

We employ now AMR to increase the resolution in correspondence of the interface. We consider the same refinement criterion \eqref{eq:eta_K} and the same thresholds for \(\eta_{K}\) adopted in Section \ref{ssec:RT} and we allow up two local refinements, so as to obtain \(h = \frac{1}{640}\) and a maximum resolution which would correspond to a \(1280 \times 2560\) uniform grid. One can notice that the resolution is enhanced close to the interface between the two fluids and a comparison of the shape of the bubble with that obtained with the uniform grid shows that we have reached grid independence (Figure \ref{fig:rising_bubble_case2_adaptive}). The final grid consists of \(274427\) elements. Figure \ref{fig:rising_bubble_case2_adaptive_comparison} reports a comparison for the quantities of interest between the fixed grid simulation, the adaptive one and the reference results. The profiles confirm that we have reached grid independence, since only the degree of circularity slightly differs between the two simulations, whereas the profiles of the vertical coordinates of both velocity and position of the center of mass are visually indistinguishable. 

\begin{figure}[!h]
	\begin{subfigure}{0.5\textwidth}
		\centering
		\includegraphics[width=0.8\textwidth]{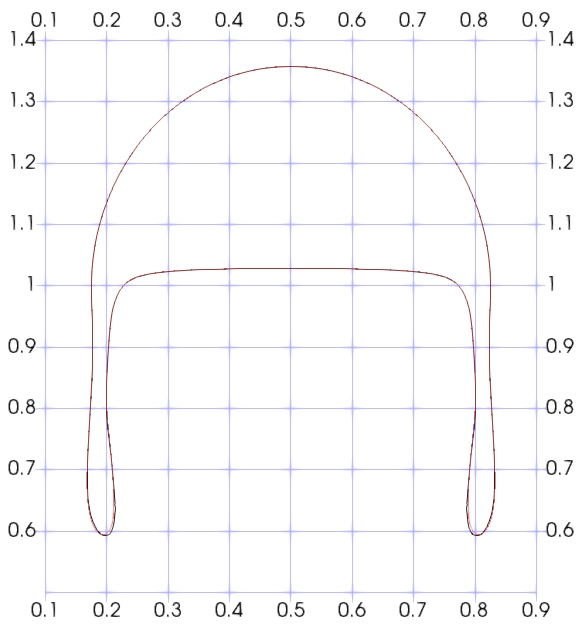}
	\end{subfigure}
	\begin{subfigure}{0.5\textwidth}
		\centering
		\includegraphics[width=0.75\textwidth]{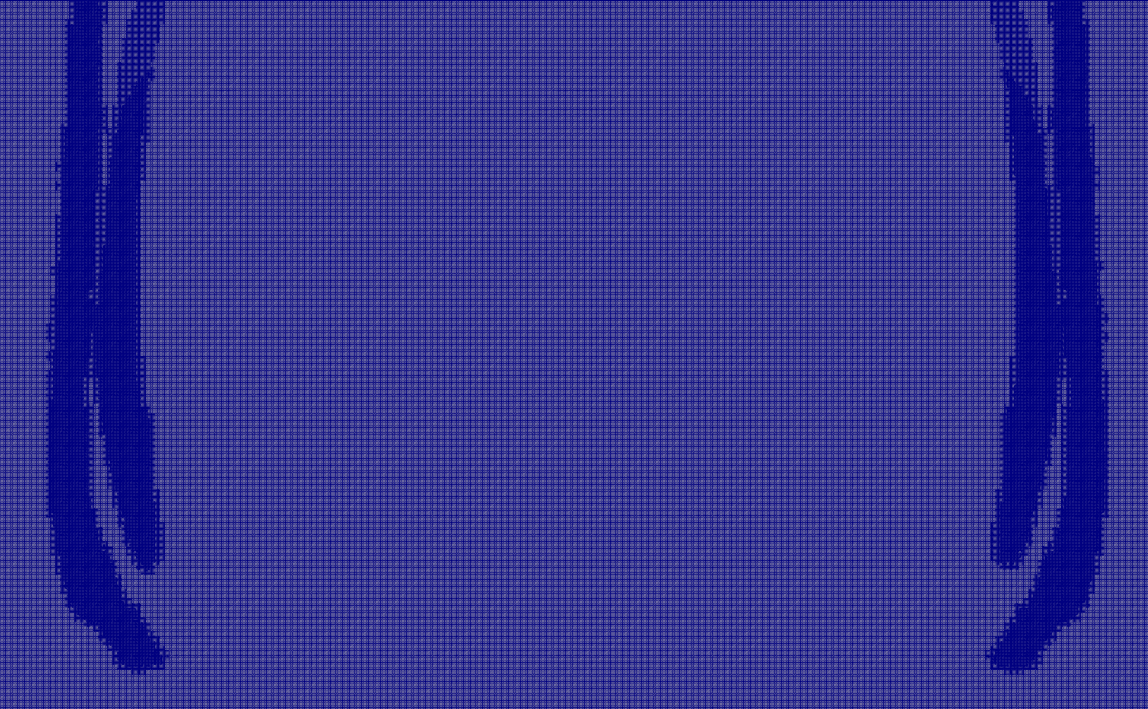}
	\end{subfigure}
	\caption{Rising bubble benchmark, configuration 2 with adaptive mesh refinement. Left: shape of bubble at $t = T_{f} = 4.2$. The black line reports the interface obtained with the adaptive mesh, while the red line shows the interface obtained with the fixed grid. Right: computational grid at $t = T_{f} = 4.2$ (close-up to thin filamentary regions).}
	\label{fig:rising_bubble_case2_adaptive}
\end{figure} 

\begin{figure}[!h]
	\begin{subfigure}{0.5\textwidth}
		\centering
		\includegraphics[width=0.9\textwidth]{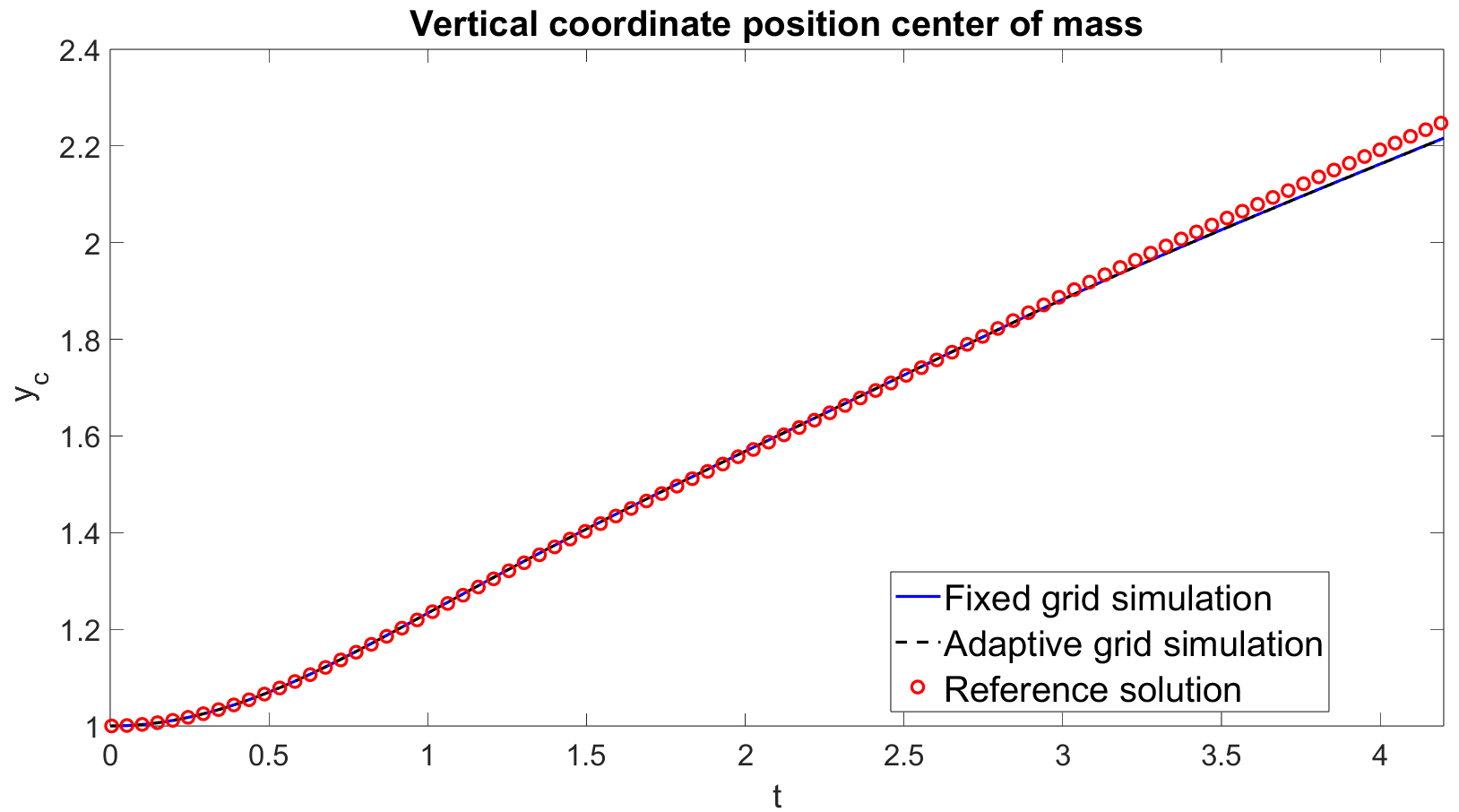} a)
	\end{subfigure}
	\begin{subfigure}{0.5\textwidth}
		\centering
		\includegraphics[width=0.9\textwidth]{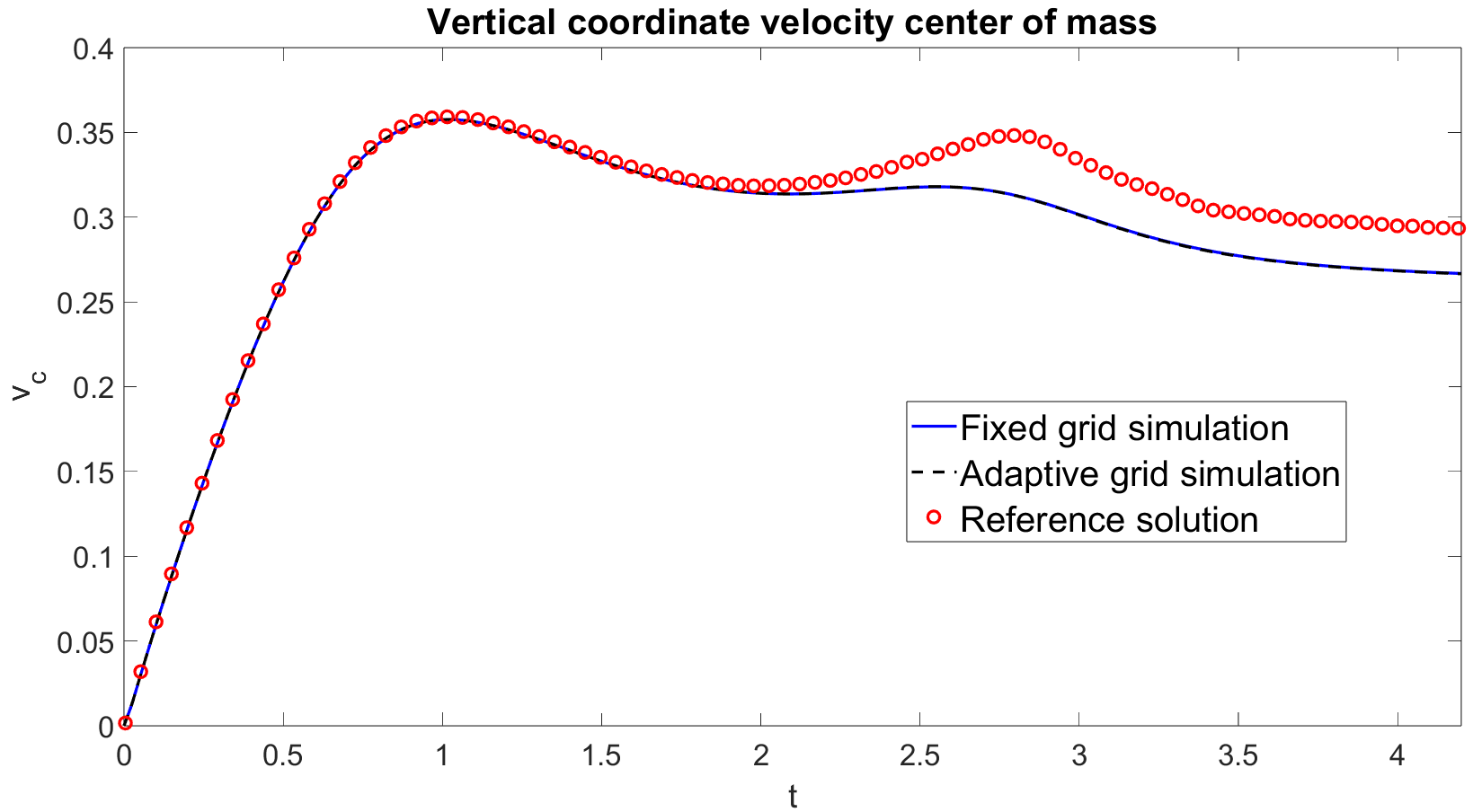} b)
	\end{subfigure}
	\begin{subfigure}{0.5\textwidth}
		\centering
		\includegraphics[width=0.9\textwidth]{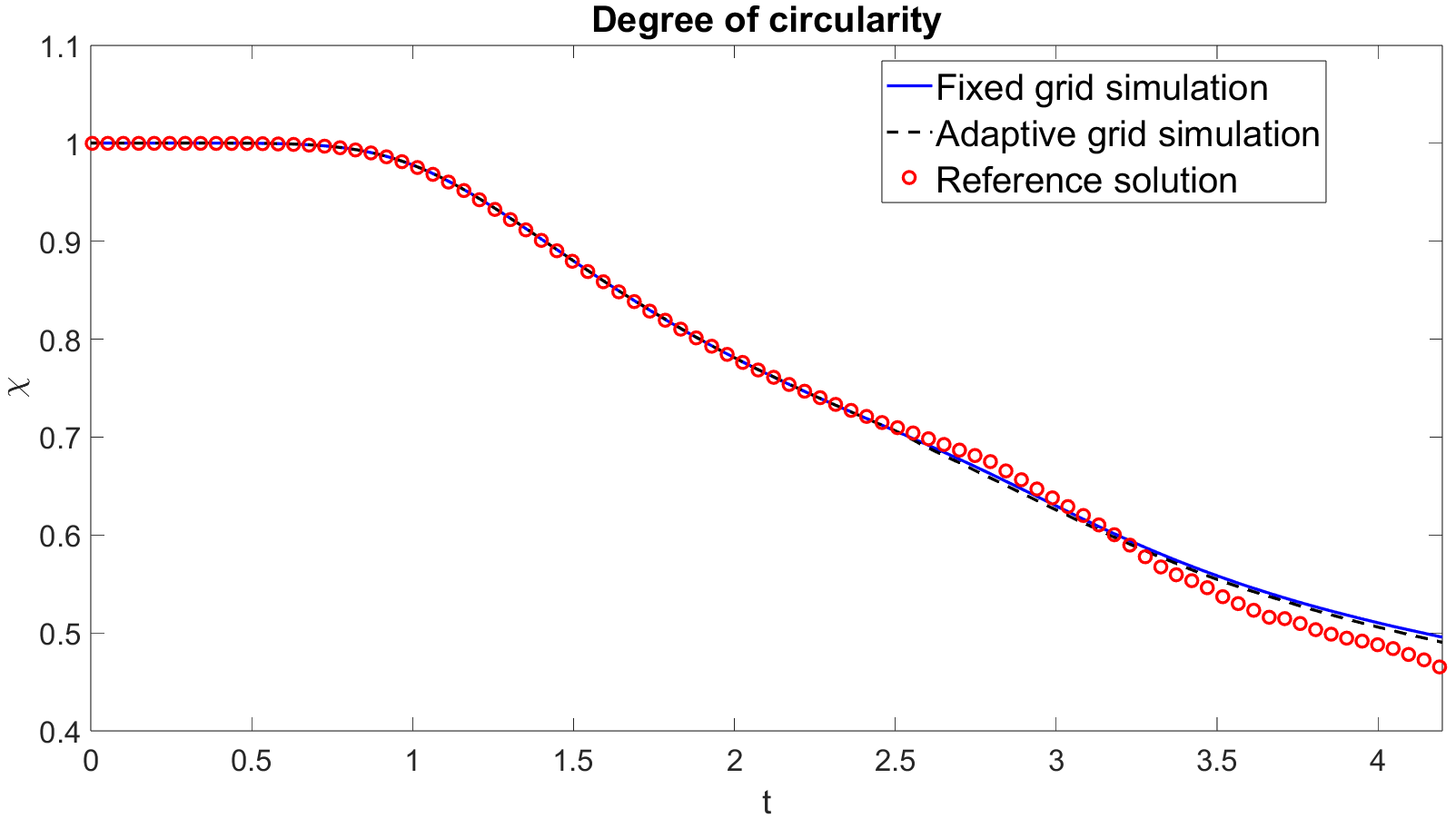} c)
	\end{subfigure}
	\caption{Rising bubble benchmark, configuration 2 with adaptive mesh refinement, a) vertical coordinate position of the center of mass, b) vertical coordinate velocity of the center of mass, c) degree of circularity. The blue lines denote the results obtained with the uniform grid, the black dashed lines report the results obtained with the adaptive grid, while the red dots show the reference results from \cite{hysing:2009}. Reference data for position and velocity have been rescaled by $L_{ref}$ and $U_{ref}$, respectively, for the sake of comparison.}
	\label{fig:rising_bubble_case2_adaptive_comparison}
\end{figure}

A significant difference in the development of the thin filamentary regions depends on the modelling of the viscosity coefficient \(\mu\), as pointed out in \cite{vrehovr:2018} for diffuse interface models. A popular alternative to the linear interpolation model defined in \eqref{eq:mu_linear} is the so-called harmonic interpolation, defined as
\begin{equation}\label{eq:harmonic_viscosity}
	\frac{1}{\mu} = H_{\varepsilon}(\varphi) + \frac{\mu_{1}}{\mu_{2}}\left(1 - H_{\varepsilon}(\varphi)\right).	
\end{equation}	
This choice yields results which are more similar to Group 1 in \cite{hysing:2009}, where a break-up occurs (Figure \ref{fig:rising_bubble_case2_shape_harmonic}). For what concerns the quantities of interest, we notice from Figure \ref{fig:rising_bubble_case2_harmonic} that, since the thin elongated filaments break themselves, the degree of circularity is higher. Moreover, both the second peak of the rising velocity and the final position of the center of mass are significantly higher and closer to the reference results. The following analysis further confirms how challenging is defining a reference benchmark solution when the bubble undergoes large deformations.

\begin{figure}[!h]
	\begin{subfigure}{0.5\textwidth}
		\centering
		\includegraphics[width=0.8\textwidth]{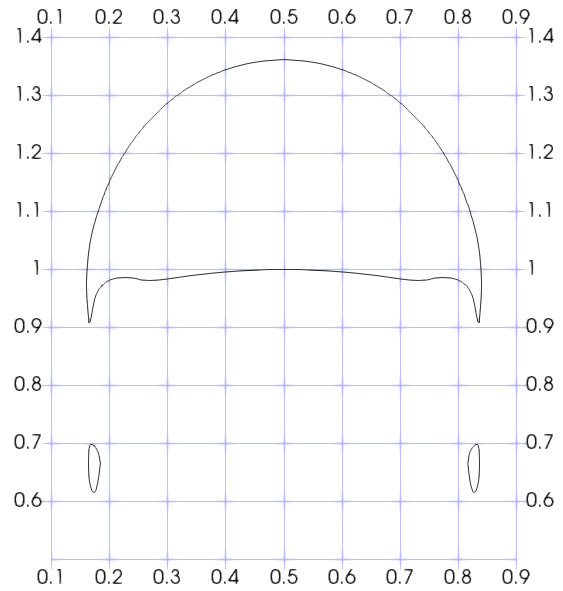}
	\end{subfigure}
	\begin{subfigure}{0.5\textwidth}
		\centering
		\includegraphics[width=0.75\textwidth]{./figures/Rising_Bubble/Case_2/shape_ref.png}
	\end{subfigure}
	\caption{Rising bubble benchmark, configuration 2, shape of bubble at $t = T_{f} = 4.2$. Left: numerical simulation using the harmonic interpolation \eqref{eq:harmonic_viscosity} for the viscosity. Right: image from \cite{hysing:2009}. Bounds have been rescaled by $L_{ref}$ for the sake of comparison with the reference results.}
	\label{fig:rising_bubble_case2_shape_harmonic}
\end{figure}

\begin{figure}[!h]
	\begin{subfigure}{0.5\textwidth}
		\centering
		\includegraphics[width=0.9\textwidth]{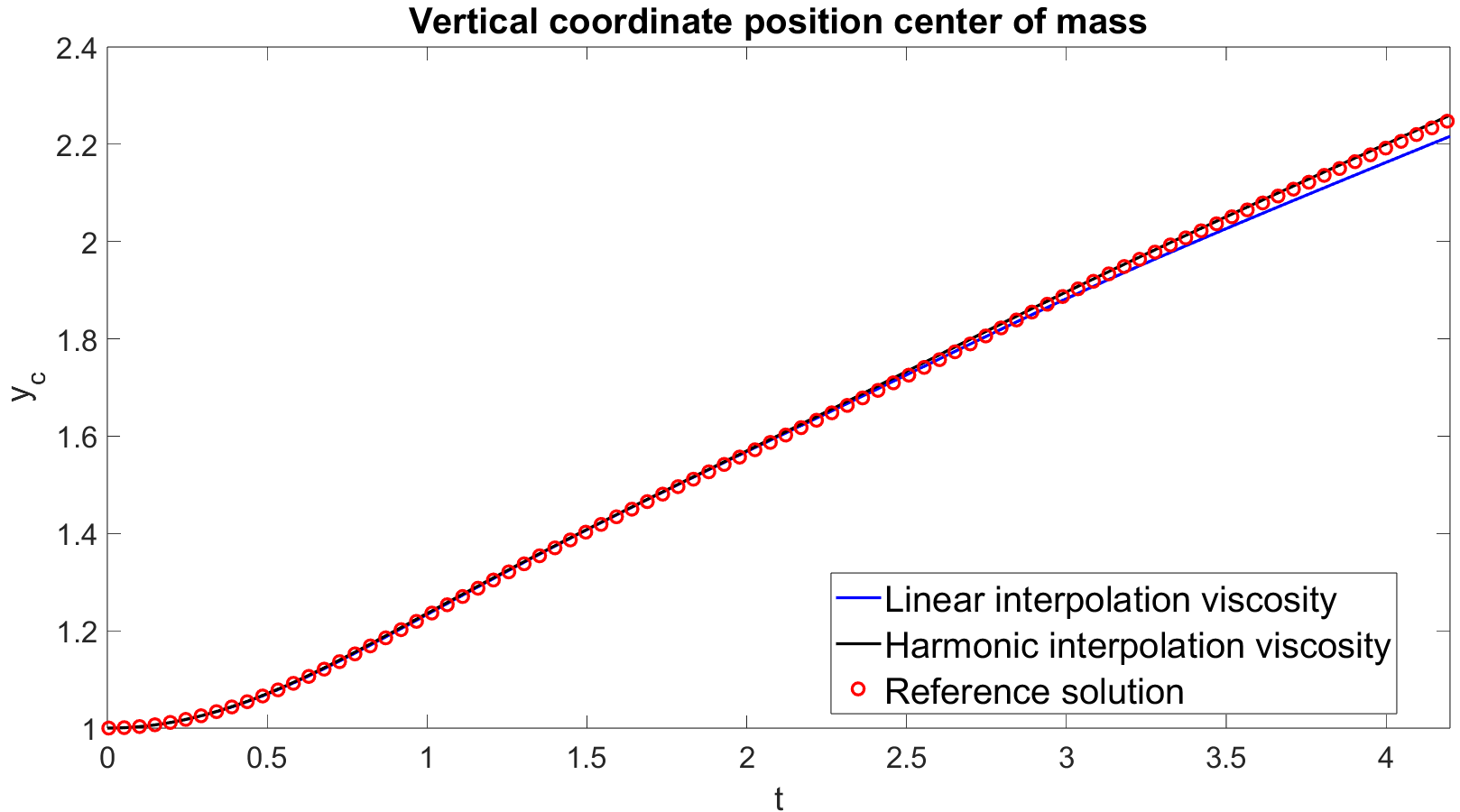} a)
	\end{subfigure}
	\begin{subfigure}{0.5\textwidth}
		\centering
		\includegraphics[width=0.9\textwidth]{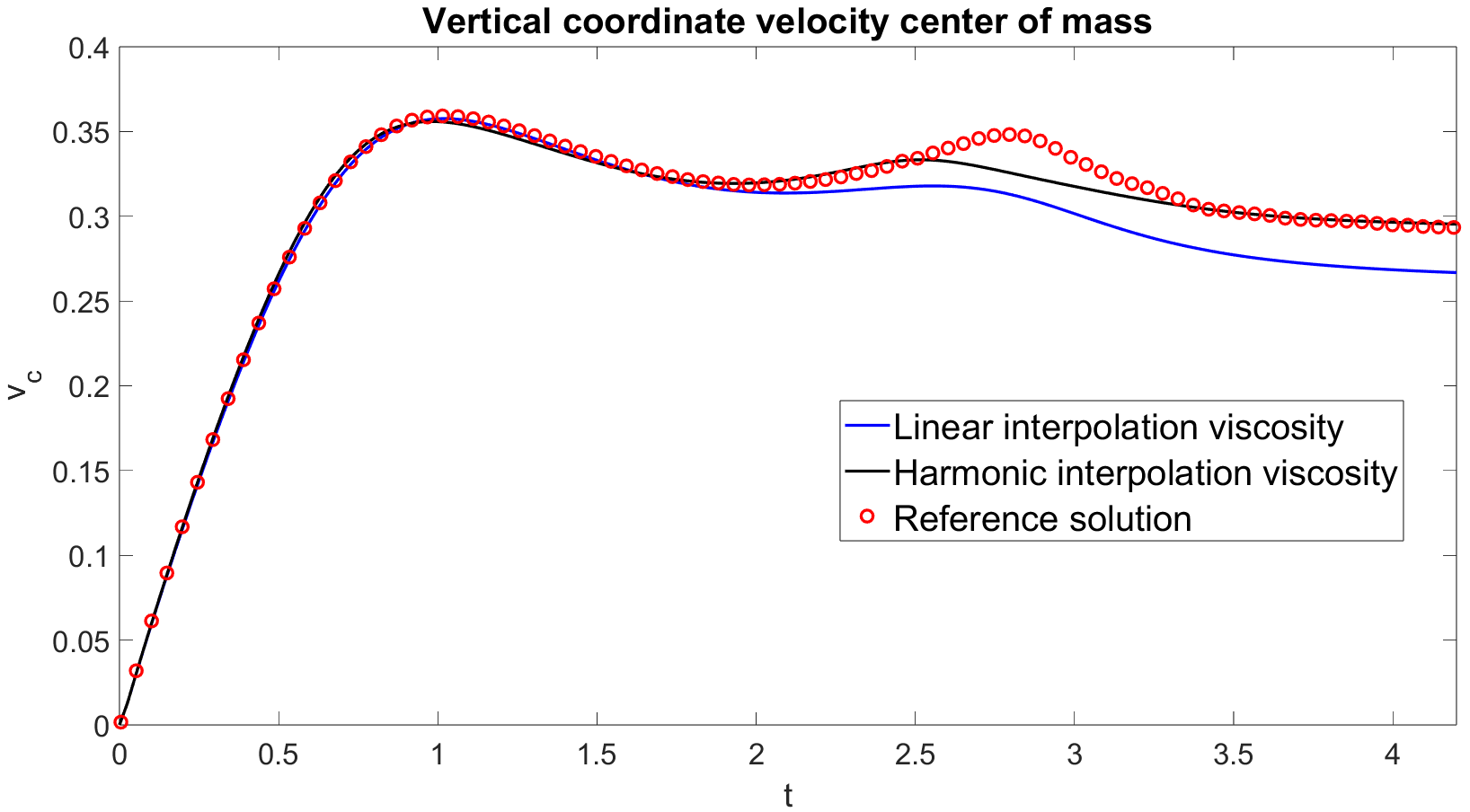} b)
	\end{subfigure}
	\begin{subfigure}{0.5\textwidth}
		\centering
		\includegraphics[width=0.9\textwidth]{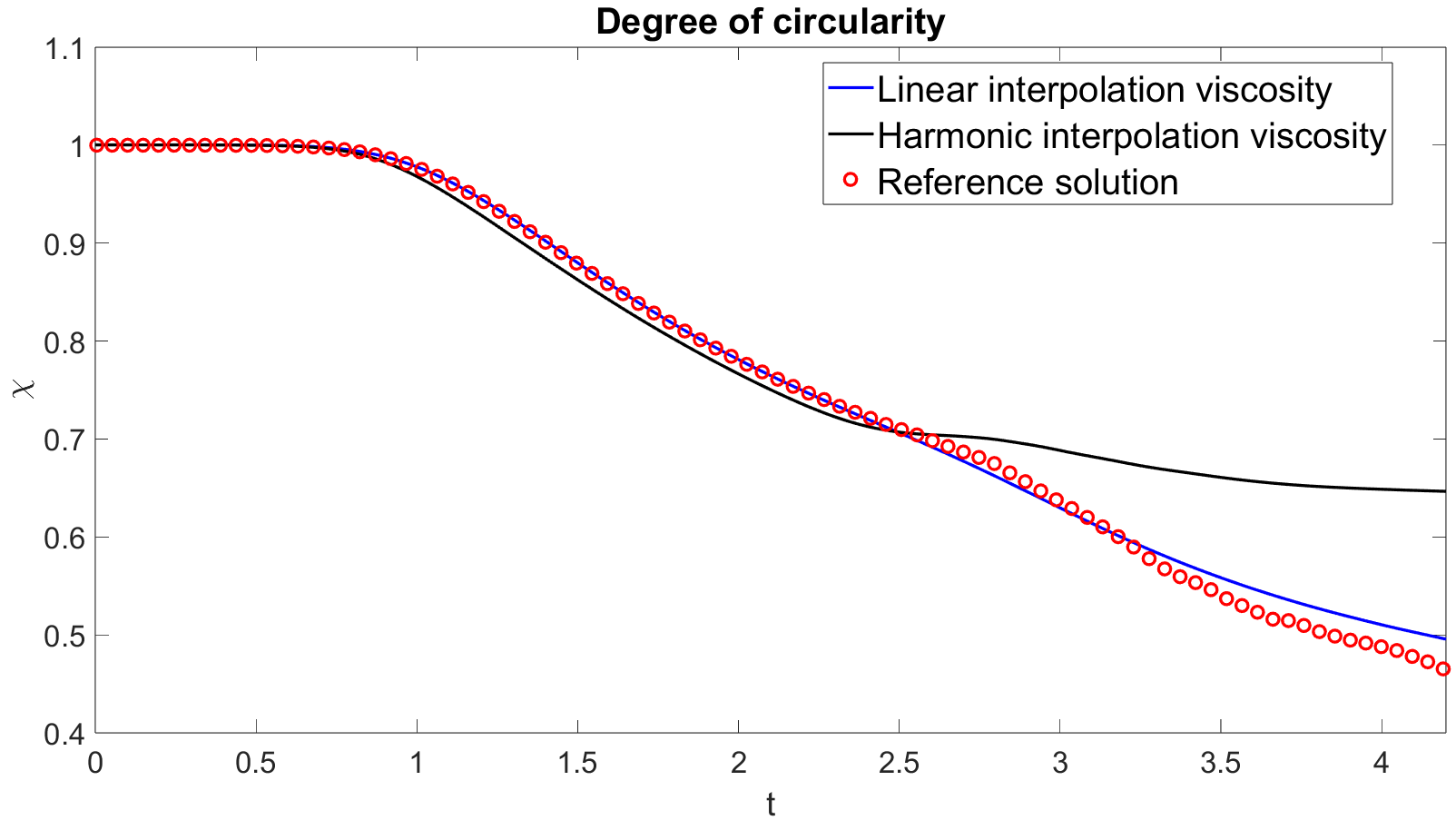} c)
	\end{subfigure}
	\caption{Rising bubble benchmark, configuration 2, comparison between linear interpolation and harmonic interpolation for the viscosity: a) vertical coordinate of the position of the center of mass, b) vertical coordinate of the velocity of the center of mass, c) degree of circularity. The blue lines denote the results with the linear interpolation, the black lines report the results with the harmonic interpolation, while the red dots show the reference results from \cite{hysing:2009}. Reference data for position and velocity have been rescaled by $L_{ref}$ and $U_{ref}$, respectively, for the sake of comparison.}
	\label{fig:rising_bubble_case2_harmonic}
\end{figure}

\subsection{Falling bubble test case}
\label{ssec:falling_bubble}

Finally, we consider the test case of a falling bubble in order to further analyze the impact of the choice of the mixture viscosity. The present configuration is a modification of that presented in \cite{grooss:2006}. The domain is \(\Omega = \SI[parse-numbers=false]{\left(0,1\right) \times \left(0, \frac{3}{2}\right)}{\meter\squared}\). The bubble is initially an ellipse  centered in \(\left(x_{0},y_{0}\right) = \SI[parse-numbers=false]{\left(0.6, 1.2\right)}{\meter}\) with semi-axes \(a = \SI{0.2}{\meter}\) and \(b = \SI[parse-numbers=false]{\frac{\sqrt{2}}{10}}{\meter}\). The density of the heavier fluid is \(\rho_{1} = \SI{1000}{\kilogram\per\meter\cubed}\), while its viscosity is \(\mu_{1} = \SI{4}{\kilogram\per\meter\per\second}\). Following \cite{grooss:2006}, we take \(L_{ref} = a = \SI{0.2}{\meter}, U_{ref} = \SI{0.5}{\meter\per\second}, \frac{\rho_{2}}{\rho_{1}} = 10^{-2}\), and \(\frac{\mu_{2}}{\mu_{1}} = \frac{1}{8}\). Hence, the computational domain is \(\Omega = \left(0, 5\right) \times \left(0, 7.5\right)\), with \(Re = 25\) and \(Fr \approx 0.35\). We take \(M = 3.3 \times 10^{-4}\), corresponding to \(c \approx \SI{1500}{\meter\per\second}\). The computational grid is composed by \(160 \times 240\) elements and we set \(\varepsilon = h = 0.03125, \Delta\tau = 0.05h, u_{c} = 0.05u_{max}\), and \(\beta = 0.5\). Surface tension is not included. Hence, the bubble will tend to break itself and we aim to show the different evolution of the interface according to the choice of the mixture viscosity. No-slip boundary conditions are imposed on top and bottom boundaries, together with homogeneous Neumann boundary conditions for the pressure. Periodic boundary conditions are prescribed in the horizontal direction. The initial velocity field is zero, while the initial level set function is described by the following relation:
\begin{equation}
	\phi(0) = \frac{1}{1 + \exp\left(\frac{\sqrt{\frac{b^{2}}{a^{2}}\left(x - \frac{x_{0}}{L_{ref}}\right)^{2}\left(y - \frac{y_{0}}{L_{ref}}\right)^{2}} - \frac{b}{L_{ref}}}{\varepsilon}\right)}.
\end{equation}
We take as final time \(T_{fin} = 0.75\), with \(\Delta t = 2.5 \times 10^{-3}\). yielding an acoustic Courant number \(C \approx 485\). One can easily notice that, at \(t = \frac{3}{5}T_{f} = 0.45\), the interface obtained with the harmonic interpolation of the viscosity \eqref{eq:harmonic_viscosity} and that obtained using the linear interpolation of the viscosity \eqref{eq:mu_linear} are visually indistinguishable and exhibit some oscillations (Figure \ref{fig:falling_bubble}, left). Finally, at \(t = T_{f}\), the two interfaces are significantly different (Figure \ref{fig:falling_bubble}, right). The interface obtained using \eqref{eq:harmonic_viscosity} shows evident oscillations and numerous small scale structures, while the interface tracked employing \eqref{eq:mu_linear} is deforming itself with reduced oscillations and less small scale structures. The presence of irregular values using the harmonic interpolation of the viscosity is confirmed by the maximum Courant number \(C_{u}\), which is \(C_{u} \approx 3.6\) using \eqref{eq:harmonic_viscosity} and \(C_{u} \approx 1\) using \eqref{eq:mu_linear}. These results further validate the considerations depicted in Section \ref{ssec:rising_bubble}: the use of the harmonic interpolation of the viscosity \eqref{eq:harmonic_viscosity} tends to enhance the break-up effect and the presence of small scale structures, while the linear interpolation of the viscosity \eqref{eq:mu_linear} is able to keep a more regular shape of the interface for longer time. Hence, the modelling of the mixture viscosity plays a key role for the description of physical phenomena in which atomization and break-up occur, like those presented e.g. in \cite{lebas:2009, vallet:1999}, or when small scale structures are explicitly considered as in \cite{loison:2023}.  

\begin{figure}[!h]
	\begin{subfigure}{0.5\textwidth}
		\centering
		\includegraphics[width=0.9\textwidth]{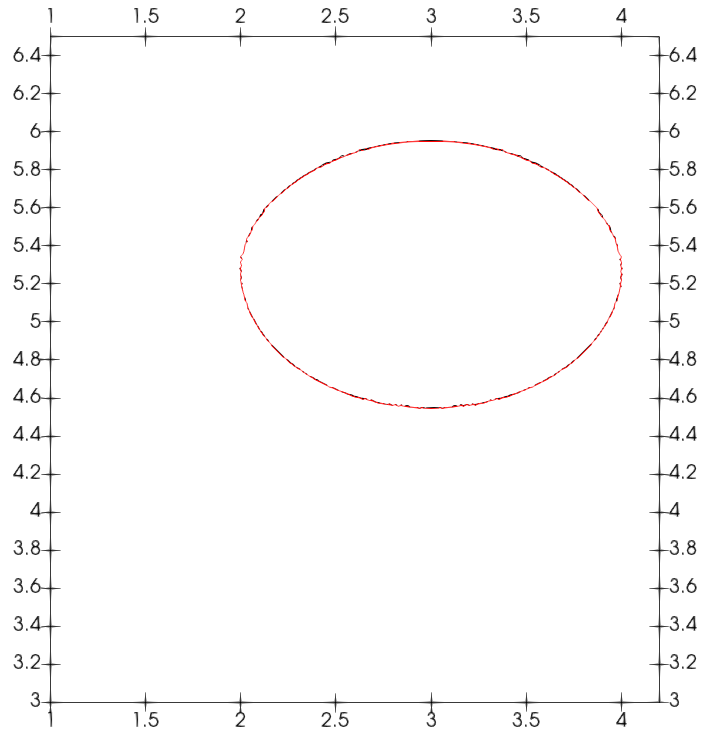}
	\end{subfigure}
	\begin{subfigure}{0.5\textwidth}
		\centering
		\includegraphics[width=0.9\textwidth]{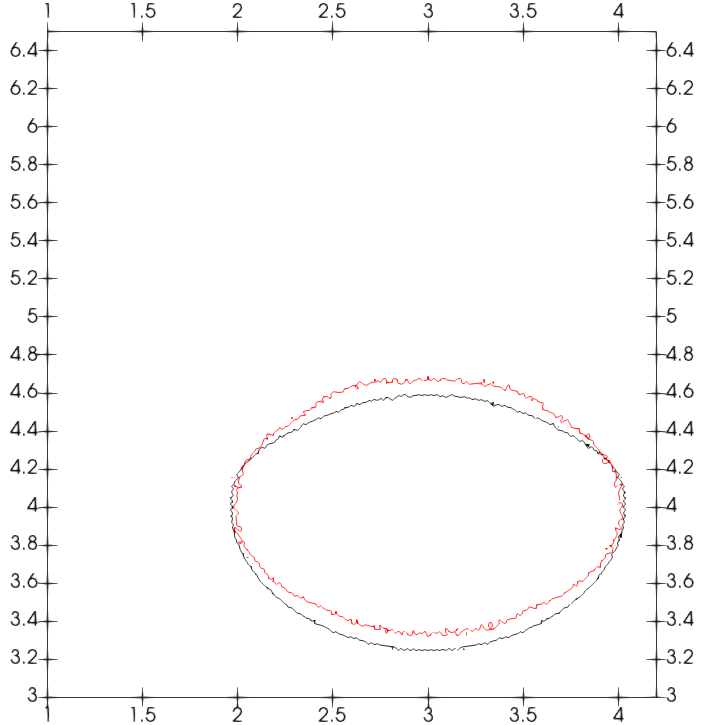}
	\end{subfigure}		
	\caption{Falling bubble test case with a $160 \times 240$ mesh. Left: shapes of the bubble at $t = \frac{3}{5}T_{f} = 0.45$. Right: shapes of the bubble at $t = T_{f} = 0.75$. The black lines report the solutions obtained using the linear interpolation for the viscosity \eqref{eq:mu_linear}, while the red lines show the results obtained using the harmonic interpolation for the viscosity \eqref{eq:harmonic_viscosity}.}
	\label{fig:falling_bubble}
\end{figure}

%%%%%%%%%%%%%%%%%%%%%%%%%%%%%%%%%%%% Conclusions %%%%%%%%%%%%%%%%%%%%%%%%%%%%%%
\section{Conclusions}
\label{sec:conclu} \indent

Building on the experience of \cite{orlando:2023b, orlando:2022}, we have proposed an implicit Discontinuous Galerkin discretization for incompressible two-phase flows. While discretizations of incompressible two-phase flows equations have been proposed in many other papers, we have presented here an approach based on an artificial compressibility formulation in order to avoid some well known issues of projection methods. The time discretization is obtained by a projection method based on the L-stable TR-BDF2 method. The implementation has been carried out in the framework of the numerical library \texttt{deal.II}, whose mesh adaptation capabilities have been exploited to increase the resolution in correspondence of the interface between the two fluids. The effectiveness of the proposed approach has been shown in a number of classical benchmarks. In particular, the influence of some possible choices for the mixture viscosity when the interface undergoes large deformations or when break-up occurs has been established, following and extending an analysis previously carried out for diffuse interface models. In future work, we aim to exploit the possibility of considering well resolved interfaces for an analysis on the evolution equations of interfacial quantities, following the recent contributions \cite{orlando:2023c, orlando:2024a}, as well as to develop an extension of analogous numerical approaches to fully compressible two-phase flows.  

\section*{Acknowledgements}

G.O. is part of the INdAM-GNCS National Research Group. The author would like to thank L. Bonaventura, P. Barbante, and E. Capuano for several useful discussions on related topics. The author gratefully acknowledges N. Parolini for providing the original data of the rising bubble test case discussed in Section \ref{sec:tests}. The simulations have been partly run at CINECA thanks to the computational resources made available through the NUMNETF-HP10C06Y02 ISCRA-C project. We acknowledge the CINECA award, for the availability of high-performance computing resources and support. This work has been partially supported by the ESCAPE-2 project, European Union’s Horizon 2020 Research and Innovation Programme (Grant Agreement No. 800897).
 
\bibliographystyle{plain}
\bibliography{DG_TRBDF2_TwoPhase}

\end{document}